\documentclass[12pt]{article}
\usepackage{amssymb}

\usepackage{amsmath}
\usepackage{latexsym}
\usepackage{amsmath, amssymb}
\usepackage{graphicx}
\usepackage{color, soul}

\newtheorem{theorem}{Theorem}
\newtheorem{lemma}{Lemma}
\newtheorem{corollary}{Corollary}

\newcommand{\mle}{\widehat{\theta}}
\newcommand{\CI}{\mbox{CI}}

\newcommand{\that}{\widehat{\theta}}

\newcommand{\half}{\frac{1}{2}}

\begin{document}

\title{Epistemic confidence, the Dutch Book and relevant subsets}
\author{{\large Yudi Pawitan, Hangbin Lee and Youngjo Lee} \\
Department of Medical Epidemiology and Biostatistics\\
Karolinska Institutet, Sweden\\
and\\
Department of Statistics\\
Seoul National University, South Korea}
\maketitle

\begin{abstract}
We use a logical device called the Dutch Book to establish epistemic confidence, defined as the sense of confidence \emph{in an observed} confidence interval. This epistemic property is unavailable -- or even denied -- in orthodox frequentist inference. In financial markets, including the betting market, the Dutch Book is also known as arbitrage or risk-free profitable transaction. A numerical confidence is deemed epistemic if its use as a betting price is protected from the Dutch Book by an external agent. Theoretically, to construct the Dutch Book, the agent must exploit unused information available in any relevant subset. Pawitan and Lee (2021) showed that confidence is an extended likelihood, and the likelihood principle states that the likelihood contains all the information in the data, hence leaving no relevant subset. Intuitively, this implies that confidence associated with the full likelihood is protected from the Dutch Book, and hence is epistemic. Our aim is to provide the theoretical support for this intuitive notion.
\end{abstract}

\section{Introduction}

Given data $Y=y$ -- of arbitrary size or complexity -- generated from a model
$p_\theta(y)$ indexed with scalar parameter of interest $\theta$, a confidence interval $\mbox{CI}(y)$ is computed with coverage probability
\begin{equation*}
P_\theta(\theta \in \mbox{CI(Y)}) = \gamma.
\end{equation*}
We are interested in the \emph{epistemic confidence}, defined as the sense of confidence \emph{in the observed CI(y).} (For simplicity, we shall often drop the
explicit dependence on $y$ from the CI.) Arguably, this is what we want from a CI, but the orthodox frequentist view is emphatic that the probability $\gamma$ does not apply to the observed interval $\CI(y)$ but to the procedure. There are well-known examples justifying this position; see Example 1 below. In the confidence interval theory, the coverage probability is called the confidence level. So, in the frequentist theory, `confidence' has actually no separate meaning from probability; in particular, it has no epistemic property. Schweder and Hjort (2016) and Schweder (2018)
have been strong proponents of interpreting confidence as `epistemic probability.' However, their view is not commonly accepted. Traditionally, only the Bayesians have no problem in stating that their subjective probability is epistemic. How do they achieve that? Is there a way to make a non-Bayesian confidence epistemic?

Frequentists interpret probability as either a long-term frequency or a propensity of the generating mechanism, such as a coin toss or a confidence interval procedure. So, for them, unique events, such as the next toss or the true status of an observed CI, do not have a probability. On the other hand, Bayesians can attach their subjective probability to such unique events. This interpretation is made possible using a logical device called the Dutch Book. As classically proposed by Frank Ramsey (1926) and Bruno de Finetti (1931), one's subjective probability of an event $E$ is defined as the personal betting price one puts on the event. Though subjective, the price is not arbitrary, but it follows a normative rational consideration; it is a price such that no external agent can construct a Dutch Book against them, i.e., make a risk-free profit. In other words, it is irrational to make a bet that is guaranteed to lose. The Dutch Book is also known as arbitrage or free lunch. In the classical Dutch Book argument, the bet is made between two individuals. 

Likewise, here \emph{we define confidence to be epistemic if it is protected from the Dutch Book,} but crucially we assume that there is a betting market of a crowd of independent and intelligent players. In this
market, bets are like a commodity with supply and demand from among the players. Assuming a perfect market condition -- for instance, full competition, perfect information and no transaction cost -- in accordance with the Arrow-Debreu theorem (Arrow and Debreu, 1954), there is an equilibrium price at which there is
balance between supply and demand. Intuitively, if you are a seller and set your price too low, many would want to buy from you, thereby creating demand and increasing the price. Whereas if you set your price too high, nobody would want to buy from you, thus reducing demand and pressuring the price down. For the betting market in particular, the fundamental theorem of asset pricing (Ross, 1976) states that, assuming an objective probability model, there is no arbitrage if the price is determined by the objective probability. `Perfect
information' here means all players have access to the generic data $y$ and
the sampling model $p_\theta(y)$. (If there is no objective probability
model, the market, as evidenced by actual betting markets, can still have an agreed interpersonal price at any given time, though not a
theoretically determined price.)

To illustrate the role of the betting market in the Dutch Book argument, suppose you and I are betting on the 2024 US presidential election. Suppose the betting market is currently giving the price of 0.25 for Donald Trump to win (this means you pay \$0.25 to get \$1 back if Trump wins, including your original \$0.25). Suppose, for whatever reasons, you believe Trump will lose and hence set the probability of him winning at 0.1. Then I would construct a Dutch Book by `buying from you' at \$0.1 and immediately `selling it in the market' at \$0.25, thus making a risk-free profit. `Buying from you' means treating you like a bookie: paying you \$0.1 and getting \$1 back if Trump wins. While `selling in the market' for me means betting
against the event, so I behave like a bookie: people pay me \$0.25 so that they can get \$1 back if Trump wins, but I keep the \$0.25 if Trump loses. So, overall I would make \$0.15 risk-free, i.e., regardless of whether Trump wins or loses. Note that this is not just a thought experiment -- you can do all this buying and selling of bets in the online betting market.

It is worth emphasizing the difference between our setup and the classical Dutch Book argument used to establish the subjective Bayesian probability. In the latter, because it does not presume the betting market, bets are made only between you and me. To avoid the Dutch Book, you have to make your bets internally consistent by following (additive) probability laws. However, \emph{even if your bets are internally consistent (or coherent),} if your prices do not match the market prices, I can make a risk-free profit by playing between you and the market; see Example 1. So, the presence of the market creates a stronger requirement for epistemic probability. We shall avoid the terms 'subjective' and 'objective'; one might consider `epistemic' to be subjective since it refers to a personal decision-making based on a unique event, but the market consideration makes it impersonal.

The present issue is when the confidence, as measured by the coverage probability, applies to the observed interval. One way to judge this is whether you are willing to bet on the true status of the CI using the confidence level as your personal price. Normatively, this should be the case if you know there is no better price. Intuitively, this is when you're sure that you have used all the available information in the data, so nobody can exploit you, i.e., construct a Dutch Book against you. Theoretically, to construct the Dutch Book, an external agent must exploit unused information in the form of a relevant subset, conditional on whether he can get a different coverage probability.

Pawitan and Lee (2021) showed that the confidence is an extended likelihood (Lee et al., 2017). The classical likelihood principle (Birnbaum, 1962) and its extended version (Bj{\o}rnstad, 1996) state that the likelihood contains all the information in the data. Intuitively, this implies that the likelihood leaves no relevant subset, and is thus protected from
the Dutch Book. In other words, we can attach the degree of confidence to the observed CI, i.e., confidence is epistemic, provided it is associated with the full likelihood. Our aim is to establish the theoretical justification for this intuitive notion.

To summarize briefly and highlight the plan of the paper, we describe three key concepts: relevant subset, confidence and ancillary statistic. We prove the main theorem that there are no relevant subsets if confidence is associated with the full likelihood. This condition is satisfied if the confidence is based on a sufficient statistic. When there is no sufficient statistic, but there is a maximal ancillary statistic, then this ancillary defines relevant subsets and there are no further relevant subsets.

\section{Main theory}

\subsection{Relevant subsets}

The idea of relevant subset appeared in Fisher's writings on the nature of probability (Fisher, 1958). He considered probability meaningful provided there is no relevant subset. However, he treated the condition as an axiom -- appealing to our intuition as to what `meaningful' means -- and did not establish what conditions are needed to guarantee no relevant
subset. To avoid unnecessary philosophical discussions, we're limiting the idea of relevant subset to the confidence interval procedures but not to the probability concept in general.

Intuitively, we could use the coverage probability $\gamma$ as a betting price if there is no better price given the data at hand. So the question is, are there any features of the data that can be used to improve the price? Mathematically,
these `features' are some observed statistics that can be used to help predict the true coverage status. Given an arbitrary statistic $S(y)$, the conditional coverage probability $P_\theta(\theta \in \mbox{CI}|S(y))$ will in general be biased, i.e., different from the marginal coverage. However, the bias as a function of the unknown $\theta$ is generally not going to be consistently in one direction. For example, trivially, if we use the full data $S(y) = y$ itself, i.e., fixing the data as observed, then the conditional coverage is either one or zero depending on the true status of the CI, hence completely non-informative. In terms of betting, this means we cannot exploit an arbitrary feature of the data as a way to construct a Dutch Book against someone who sets the price at $\gamma$. The betting motivation also appeared in Buehler (1959) and Robinson (1979), though they only assumed two people betting each other repeatedly, but not the existence of the betting market. As we shall discuss in Example~1 and after Theorem~1, this has a significant impact on the interpretation of epistemic confidence.

A statistic $R(y)$ is defined to be relevant (cf.\ Buehler, 1959) if the conditional coverage is \emph{non-trivially and consistently biased in one direction.} That is, for a positive bias,
there is $\epsilon> 0$ free of $\theta$, such that
\begin{eqnarray}
P_\theta(\theta \in \mbox{CI}(Y)|R(y)) \ge \gamma+ \epsilon
\quad \text{for all } \theta.
\label{eq:bias}
\end{eqnarray}
Now, potentially the feature $R(y)$ \emph{can be used} to construct a Dutch Book: Suppose you and I are betting, and I notice that the event $[R(y)]$
occurs. If you set the price at $\gamma$, then I would buy the bet from you and then sell it in the betting market at $\gamma + \epsilon$. So I make
a risk-free profit of $\epsilon$. (We have assumed that the market contains
intelligent players, so they would also have noticed the relevant statistic
and set the price accordingly.) Similarly, for the negative bias, the
relevant $R(y)$ has the property
\begin{equation}
P_\theta(\theta \in \mbox{CI}(Y)|R(y)) \le \gamma- \epsilon\ \mbox{for all} \ \theta.
\label{eq:negative bias}
\end{equation}

Technically, $R(y)$ induces subsets of the sample space, known as the `relevant subset'; for convenience, we use the terms 'relevant statistic' and 'relevant subset' interchangeably.  So, if there is a relevant subset, the confidence level $\gamma$ is not epistemic. Conversely, if there are no relevant subsets, the betting price determined by the confidence level is protected from the Dutch Book. So, \emph{mathematically, we establish epistemic confidence by showing that it corresponds to a coverage probability that has no relevant subsets.} \vspace{.1in}

\noindent \textbf{Example 1.} Let $y\equiv (y_1,y_2)$ be an iid sample from a uniform distribution on $\{\theta-1,\theta,\theta+1\}$, where the parameter $\theta$ is an integer. Let $y_{(1)}$ and $y_{(2)}$ be
the minimum and maximum values of $y_1$ and $y_2$. We can show that the
confidence interval $\mbox{CI}(y) \equiv[y_{(1)},y_{(2)}]$ has a coverage
probability
\begin{equation*}
P_\theta(\theta\in \mbox{CI}) = 7/9 = 0.78.
\end{equation*}
For example, on observing $y_{(1)}=3$ and $y_{(2)}=5$, the interval $[3,5]$ is formally a 78\% CI for $\theta$. But, if we ponder a bit, in this case we can
actually be sure that the true $\theta=4$. So,
the probability of 7/9 is clearly a wrong price for this interval. This is a typical example justifying the frequentist objection to attaching the coverage probability as a sense of confidence in an observed CI.

Here the range $R\equiv R(y) \equiv y_{(2)}-y_{(1)}$ is
relevant. If $R=2$ we know for sure that $\theta$ is equal to the midpoint
of the interval, so the CI will always be correct. But if $R=0$, the CI is
equal to the point $y_1$, and it falls with equal probability at the
integers $\{\theta-1,\theta,\theta+1\}$. So, for all $\theta$, we have
\begin{eqnarray*}
P_\theta(\theta\in \mbox{CI}|R=2) &=& 1 > 7/9 \\
P_\theta(\theta\in \mbox{CI}|R=1) &=& 1 > 7/9 \\
P_\theta(\theta\in \mbox{CI}|R=0) &=& 1/3 < 7/9.
\end{eqnarray*}
In the betting market, the range information will be used by the intelligent players to settle prices at these conditional probabilities. We can be sure, for example,
that if $y_1=3$ and $y_2=5$, the intelligent players will not use 7/9 as the price and will instead use 1.00. So, the information can be used to construct a Dutch Book against anyone who ignores $R$ and unwittingly uses the unconditional coverage. How do we know that there is a relevant subset in this case? Moreover, given $R$, how do we know if there is no further relevant subset? 

To contrast with the classical Ramsey-de Finetti Dutch Book argument, suppose $y_1=y_2=3$. If for whatever subjective reasons, you set the price 7/9 for $[\theta\in \CI]$, you are being internally consistent as long as you set the price 2/9 for $[\theta \not\in \CI]$, since the two numbers constitute a valid probability measure. Internal consistency means that I cannot make a \emph{risk-free} profit from you based on \emph{this single realization} of $y$. Even if I know based on the conditional coverage that 1/3 is a better price, I cannot take any advantage of you because there is no betting market. So 7/9 is a valid subjective probability.

Now consider Buehler-Robinson's setup, again assuming no betting market and supposing $y_1=y_2=3$. They would say the marginal price 7/9 is a bad idea, because there is a relevant subset giving a conditional probability 1/3. In a series of independent repeated bets, if you set the price 7/9 whenever $y_1=y_2=3$, I will be happy to `sell' you the bet and be guaranteed to win \emph{in long term.} This is the usual frequentist interpretation; \emph{in any single bet} I am not guaranteed a risk-free money. As previously described, the presence of the betting market allows me make free-money from a single realization of $y$. So, the threat of the market-based Dutch Book is more potent. The exact technical difference between Buehler-Robinson's and our setup will be discussed below after Theorem 1, where the former allows one to choose an arbitrary prior distribution.
$\Box$ \vspace{.1in}

\subsection{Confidence distribution}

It turns out that establishing a no-relevant-subset condition
relies on the concept of a confidence. Let $t\equiv T(y)$ be a
statistic for $\theta $, and define the right-side P-value function
\begin{equation}
C_{m}(\theta ;t)\equiv P_{\theta }(T\geq t).  \label{eq:uniform}
\end{equation}
Assuming that, for each $t$, it behaves formally like a proper cumulative
distribution function, $C_{m}(\theta ;t)$ is called the confidence
distribution of $\theta $. The subscript $m$ is used to indicate that it is a `marginal' confidence, as it depends on the marginal distribution of $T$.
For continuous $T$, at the true parameter, the random variable $C_{m}(\theta
;T)$ is standard uniform. For continuous $\theta $, the corresponding
confidence density is
\begin{equation}
c_{m}(\theta )\equiv c_{m}(\theta ;t)\equiv \partial C_{m}(\theta
;t)/\partial \theta.
\label{eq:ctheta}
\end{equation}
The functions $C_{m}(\theta ;t)$ and $c_{m}(\theta )$ across $\theta $ are
\emph{realized statistics}, which depend on both the data and the model, but not on the true unknown parameter $\theta _{0}$. We can view the confidence distribution simply as the collection of P-values or CIs. Suggestively, and with a slight abuse of notation, we define
\begin{equation}
C_{m}(\theta \in \mbox{CI})\equiv \int_{\textrm{CI}}c_{m}(\theta )d\theta
\label{eq:marginal-conf}
\end{equation}
to convey the `confidence of $\theta $ belonging in the CI'.

We assume a regularity condition, called R1 below, that for any $\alpha\in (0,1)$, the quantile function $q_\alpha (\theta)$ of $T$ is a strictly increasing function of $\theta$. Then the frequentist procedure based on $T$ gives
a $\gamma$-level CI defined by
\begin{equation}
\label{ci}
\CI_\gamma(T) = \left( q_{\gamma_2}^{-1}(T), q_{\gamma_1}^{-1}(T) \right)
\end{equation}
for some $\gamma_2>\gamma_1>0$ with $\gamma_2-\gamma_1=\gamma$,
to have a coverage probability
$$
P_\theta (\theta \in \CI_\gamma (T)) = P_\theta \Big[ T \in (q_{\gamma_1}(\theta), q_{\gamma_2}(\theta)) \Big] = \gamma_2 -\gamma_1 =\gamma.
$$
Here the coverage probability is a frequentist probability based on the distribution of $T$,
whereas the confidence is for the observed interval $\CI(t)$ based on the confidence density of $\theta$.
The confidence becomes
\begin{align*}
C_m (\theta \in \CI_\gamma(t); t) &= C_m(\theta=q_{\gamma_1}^{-1}(t); t) - C_m(\theta=q_{\gamma_2}^{-1}(t); t) \\
&= P_{\theta=q_{\gamma_1}^{-1}(t)}(T\geq t) - P_{\theta=q_{\gamma_2}^{-1}(t)}(T\geq t) \\
&= (1-\gamma_1) - (1-\gamma_2) = \gamma \\
&=P_\theta (\theta \in \CI_\gamma (T)) .
\end{align*}
Thus, we have the following lemma.
\begin{lemma}
Under the regularity condition R1,
\begin{equation} \label{eq:cover}
P_\theta (\theta \in \CI(T)) = C_m (\theta \in \CI(t); t).
\end{equation}
where $\CI(t)$ is the observed interval of confidence procedure $\CI(T)$ defined in \eqref{ci}.
\end{lemma}
Fisher (1950) was against the idea of interpreting the level of significance as a long-term frequency
in repeated samples from the same population. For $i=1,\cdots ,n,$
suppose $T_{i}$'s are estimates for $\theta _{i}$'s from different populations$.$ Let $X_{i}=I(\theta _{i}\in CI(T_{i}))$ such that $\gamma
=C_{m}(\theta _{i}\in CI(t_{i}))$ and let $\bar{X}=\sum X_{i}/n.$ Then, $%
\bar{X}=\gamma +O_{p}(1/\sqrt{n}),$ so that $\gamma$ can be a long-term frequency of true coverage from different populations or experiments.
\vspace{.1in}

\noindent \textbf{Example 2.} On observing $t$ from $N(\theta
,1)$, we have the confidence distribution
\begin{equation*}
C_m(\theta;t)=P_\theta(T\geq t)=1-\Phi (t-\theta ),
\end{equation*}
where $\Phi (\cdot )$ is the standard normal distribution function, with
corresponding confidence density $c_m(\theta )=\phi (t-\theta )$, i.e., the
normal density centered at $\theta =t$. In principle, we can derive any
confidence interval or P-value from this confidence density. This example
applies in most large sample situations where, under regular conditions, the
normal model is correct asymptotically. Furthermore, it illustrates clearly
the canonical relationship between confidence and coverage probability. For
instance, for a 95\% CI, we have
\begin{equation*}
C_m(\theta \in \mbox{CI(t)}) = 95\%,
\end{equation*}
reflecting the 95\% confidence that the observed CI covers the true
parameter. This confidence is associated with an exact coverage probability
\begin{equation*}
P_\theta(\theta\in \mbox{CI(T)})=0.95,
\end{equation*}
so the confidence matches the coverage probability.$\Box$ \vspace{.1in}

Fisher (1930, 1933) called $C_m(\theta;t)$ the fiducial distribution of $\theta$, but he required $T$ to be sufficient. However, the recent definition of the confidence distribution (e.g. Schweder and Hjort,
2016, p.58) requires only $C_m(\theta;T)$ to be uniform at the true parameter, thus guaranteeing a correct coverage probability. Lemma 1 states when Fisher's fiducial probability $C_m(\theta;t)$ becomes a frequentist probability, which requires $T$ to be continuous. When $T$ is discrete, the equality is only achieved asymptotically; see Appendix A3 for an example.

However, as shown in Example 1, a correct coverage probability does not rule out relevant
subsets. This means that the current definition of the confidence distribution does not guarantee epistemic confidence. The key step is to define a confidence distribution that uses the full information. Motivated by the Bayesian formulation and Efron (1993), first define the implied prior as
\begin{equation}
c_0(\theta)\equiv c_0(\theta;t) \equiv m(t)\frac{c_m(\theta;t)}{L(\theta ;t)},
\label{eq:c0}
\end{equation}
where $m(t)$ cancels out all the terms not involving $\theta$ in ${c_m(\theta;t)}/{L(\theta ;t)}$. In this paper, the full confidence density is defined by
\begin{equation}
c_f(\theta)  \equiv c_f(\theta;y)  \propto c_0(\theta) L(\theta;y).
\end{equation}
The subscript $f$ is now used to indicate that it is associated with the full likelihood based on the whole data. When necessary for clarity, the dependence of the confidence density and the likelihood on $t$ and on the whole data $y$ will be made explicit. $c_f(\theta)$ is defined only up to a constant term to allow it to integrate to one. Obviously, if $T$ is sufficient, then $c_m(\theta) = c_f(\theta)$, but in general they are not equal. In Section 3, we show a more convenient way to  construct $c_f(\theta)$. The confidence parallel to (\ref{eq:marginal-conf}) can be denoted by $C_{f}(\cdot )$. Thus, the full confidence density looks like a Bayesian posterior. However, the implied prior is not subjectively selected, and can be improper and data-dependent.

The full confidence density $c_f(\theta)$ is used in general to compute the degree of confidence $\gamma$ to any \emph{observed} $\CI(y)$ as
$$
\gamma = \int_{CI(y)}c_f(\theta)d\theta.
$$
The CI has a coverage probability, which may or may not be equal to $\gamma$. We say that $c_f(\theta)$ has no relevant subsets, if there is no $R(y)$ such that the conditional coverage probability is biased in one direction according to (\ref{eq:bias}) or (\ref{eq:negative bias}).  

For our main theorem, we assume the following regularity conditions, the proof is given in  Appendix A1. For completeness and easy access, R1 is restated in full here.
\begin{itemize}
\item[\textbf{R1.}]  $T=T(Y)$ is a continuous scalar statistic whose
quantile function $q_{\alpha }(\theta )$, defined by
\begin{equation*}
P_{\theta }(T\leq q_{\alpha }(\theta ))=\alpha,
\end{equation*}
is strictly increasing function of $\theta $ for any $\alpha \in (0,1)$.
\item[\textbf{R2.}] $c_{0}(\theta)$ is positive and
locally integrable on the parameter space $\Theta$ such that
\begin{equation*}
\int_{J}c_{0}(\theta )d\theta <\infty ,\ \ \text{for any compact subsets}\ \
J\subseteq \Theta .
\end{equation*}

\item[\textbf{R3.}] $\log c_{0}(\theta)$ is uniformly continuous in $y$.

\item[\textbf{R4.}] The confidence interval $\mbox{CI}(y)=(b_{L}(y),b_{U}(y))$ is \textit{locally bounded}, i.e., for any compact
set $K$ in the sample space of $y$, there exist $M_{1}$ and $M_{2}$ such that
\begin{equation*}
|b_{L}(y)|\leq M_{1}\ \text{and}\ |b_{U}(y)|\leq M_{2} \quad \text{for any } y\in K
\end{equation*}
\end{itemize}

\begin{theorem}
Consider the full confidence density $c_f(\theta) \propto c_0(\theta) L(\theta;y)$,
with $c_0(\theta)$ being the implied prior defined by (\ref{eq:c0}) satisfying R2 and R3, based on $T(Y)$ that satisfies R1. Let $\gamma$ be the degree of confidence for the observed confidence interval $\CI(y)$ that satisfies R4, such that
$$
\gamma = \int_{\textrm{CI}(y)} c_f(\theta) d\theta,  \quad \text{for all } y,
$$
Then $c_f(\theta)$ has no relevant subsets.
\end{theorem}

Note we have two ways of computing the price of an observed CI: using $C_f(\theta\in \CI)$ or using $P_\theta(\theta\in \CI)$. The latter is not guaranteed to be free of relevant subsets, while the former is not guaranteed to match the coverage probability. If the two are equal, we have a confidence that corresponds to a coverage probability that has no relevant subsets, hence epistemic confidence. If $T$ is sufficient and satisfies R1, Lemma 1 implies that the frequentist CI satisfies
$$
P_{\theta }(\theta \in \mbox{CI}(Y))= C_m(\theta \in \mbox{CI}(y))= C_f(\theta \in \mbox{CI}(y)) = \gamma,\ 
\text{for all } \theta \text{ and } y
$$
Thus, we summarize the first key result in the following corollary:
\begin{corollary}
Under the regularity conditions R1-R4, if $T$ is sufficient statistic, the confidence based on $c_m(\theta;t)$ has a correct coverage probability and no relevant subset. Hence the confidence is epistemic.
\end{corollary}

In Example 2, on observing $t$ from $N(\theta ,1)$, $c_m(\theta) =
c_f(\theta) = \phi(t-\theta). $ Furthermore, we also have $L(\theta) = \phi(t-\theta)$, so the implied prior $c_0(\theta) = 1$. The coverage probabilities match the full confidence, and by Corollary 1, the confidence is epistemic.

We note that $P_{\theta }(\theta \in \CI(Y)) = C_{f}(\theta \in \CI(y))$ holds asymptotically, regardless whether $y$ is continuous or discrete. Corollary 1 specifies the conditions where it is true in finite samples.

For more generality, it is actually more convenient to prove the theorem using $c_0(\theta)\equiv c_0(\theta;y)$ an arbitrary function that satisfies R2 and R3, as long as it leads to a proper $c(\theta;y)$. In particular, it does not have to be an implied prior (\ref{eq:c0}) that depends on the statistic $T$. If $c_0(\theta)$ is a proper probability density that does not depend on $y$,  then $c_f(\theta)$ is a Bayesian posterior density, shown already by Robinson's (1979) Proposition 7.4 not to have relevant subsets. For proper priors, $\log c_0(\theta)$ is trivially uniformly continuous in $y$, so the theorem extends his result to improper and data-dependent priors.

However, there is a significant impact on the interpretation. If you use an arbitrary $c_{0}(\theta ;y)$ that is not the same as the implied prior, \emph{and there is a betting market}, your price $\gamma$ will differ from the market price. So, as illustrated in the Introduction and in Example 1, I can construct a Dutch Book against you. This means that, in this case, the theorem is meaningful only for two people betting repeatedly against each other, with gains or losses expressed in terms of expected value or long-term average. It is exactly the setting described by Buehler (1959) and Robinson (1979). Crucially, in such a setting, the presence of relevant subsets does not guarantee an external agent to make a risk-free profit from a single bet. In this sense, it does not satisfy our original definition of epistemic confidence.


Lindley (1958) showed that, assuming $T$ is sufficient, Fisher's fiducial probability -- hence the marginal confidence -- is equal to the Bayesian posterior if and only if the family $p_\theta(y)$ is transformable to a location family. However, his proof assumed $c_{0}(\theta)$ to be free of $y$. Condition R3 of the theorem allows $c_0(\theta)$ to depend on the data, so our result is not limited to the location family.

\subsection{Ancillary statistics}
The current definition of confidence (e.g. Schweder  and  Hjort,  2016,  p.58) only requires $C_m(\theta; T)$ to follow uniform distribution. However, if $T$ is not sufficient, the marginal confidence is not epistemic, because it does not use the full likelihood, so it is not guaranteed free of relevant subsets. Limiting ourselves to models with sufficient statistics to get epistemic confidence is overly restrictive, since sufficient statistics exist at arbitrary sample sizes in the full exponential family only (Pawitan, 2001, Section 4.9). Using non-sufficient statistics implies a potential loss of efficiency and epistemic property. Further progress depends on the ancillary statistic, a feature or a function of the data whose distribution is free of the unknown parameter. As reviewed by Ghosh et al.\ (2010), it is one of Fisher's great insights from the 1920s.
We first have a parallel development for the conditional confidence distribution given the ancillary $A(y)=a$:
\begin{eqnarray*}
C_c(\theta; t|a) &\equiv & P_\theta(T \ge t|a) \\
 c_c(\theta;t|a )& \equiv& \partial C_c(\theta;t|a)/\partial \theta.
\end{eqnarray*}
As for the marginal case, we have the following corollary from Lemma 1.  Condition R1 needs a little modification, where it refers to the conditional statistic $T|a$ for each $a$.
\begin{corollary}
Under the regularity condition R1,
\begin{equation} \label{eq:cover=cc}
P_\theta (\theta \in \CI|a) = C_c(\theta \in \CI; t|a).
\end{equation}
where $\CI$ is the confidence interval based on the conditional distribution of $T|a$.
\end{corollary}
Furthermore, define the implied prior as
\begin{equation}
c_0(\theta) \equiv c_0(\theta; t|a) \equiv m(t,a) \frac{c_c(\theta;t|a)}{L(\theta;t|a)},\label{eq:c0cc}
\end{equation}
where $m(t,a)$ cancels out all the terms not involving $\theta$ in ${c_c(\theta;t|a)}/{L(\theta;t|a)}$. As before, the full confidence is $c_f(\theta) \propto c_0(\theta)L(\theta;y).$

Suppose $T(y)=t$ is not sufficient but $(t,a)$ is, where $a$ is an ancillary statistic.
In this case, $a$ is called an ancillary complement, and in a qualitative sense it is a maximal ancillary, because
\begin{eqnarray}
L(\theta ;y) &=&L(\theta ;t,a)  \notag \\
&\propto& p_{\theta }(t|a)p(a)  \notag \\
&\propto &p_{\theta }(t|a)=L(\theta ;t|a).  \label{eq:cond}
\end{eqnarray}
Thus, conditioning a non-sufficient statistic by a maximal ancillary has recovered the lost information and restored the full-data likelihood.
In particular, the conditional confidence becomes the full confidence: $c_c(\theta;t|a) = c_f(\theta)$.
Note that (\ref{eq:cond}) holds for any maximal ancillary,
so if a maximal ancillary exists, then the full likelihood is automatically equal to the conditional likelihood given any maximal ancillary statistic. In its sampling theory form, when $t$ is the maximum likelihood estimator (MLE) $\hat{\theta},$ full information can be recovered from $p_{\theta }(\hat{\theta}|a),$ whose approximation has been studied by Barndorff-Nielsen (1983).

In conditional inference (Reid 1995), it is commonly stated that we condition on the ancillary to make our inference more `relevant' to the data at hand, in other words, more epistemic. But this is typically stated on an intuitive basis; the following corollary provides a mathematical justification. Since we already condition on $A(y)$, a \emph{further
relevant subset} $R(y)$ is such that the conditional probability $P_\theta(\theta\in \mbox{CI}|A(y),R(y))$ is non-trivially and consistently
biased in one direction from $P_\theta(\theta\in \mbox{CI}|A(y))$ in the same manner as (\ref{eq:bias}). As we describe following Theorem 1 above, the result holds for an arbitrary $c_0(\theta)$ that satisfies R2-R3 and leads to a valid confidence density. So it applies to $c_0(\theta)$ defined by (\ref{eq:c0cc}). Following a similar reasoning as for the previous corollary, we can state our second key result:


\begin{corollary}
If $A(y)=a$ is maximal ancillary for $T(y)$ and $\CI$ is constructed from the conditional confidence density based on $T|a$, then under R1-R4, the conditional confidence $C_c(\theta \in \CI; t|a)$ has a correct coverage probability and no further relevant subsets. Hence the conditional confidence is epistemic.
\end{corollary}
\emph{Remark:} In view of (\ref{eq:cond}), the confidence is epistemic for any choice of the maximal ancillary. Basu (1959) showed under mild conditions that maximal ancillaries exist. However, they may not be unique; this is an issue traditionally considered most problematic in conditional inference. If the maximal ancillary is not unique, then the conditional coverage probability might depend upon the choice. However, this does not affect the absence of relevant subset guaranteed by the corollary. We discuss this further in Section~\ref{sec:disc} and illustrate with an example in Appendix A4.

\section{\label{sec:examples}Examples}
Our overall theory suggests that, regardless of the existence of a sufficient statistic, we can get epistemic confidence by computing CIs based on the full confidence density $c_f(\theta) \propto c_{0}(\theta )L(\theta ;y)$. The corresponding coverage probability is either a marginal probability or a conditional probability given a maximal ancillary, depending on whether there exists a scalar sufficient statistic. The full likelihood $L(\theta ;y)$ is almost always easy to compute. However, in order to get a correct coverage,  $c_0(\theta)$ is defined by (\ref{eq:c0}) or (\ref{eq:c0cc}), which in practice can be difficult to evaluate. For example, if we use the MLE, in general it has no closed form formula, and computing the P-value, even from its approximate distribution based on Barndorff-Nielsen's (1983) formula, can be challenging. We illustrate through a series of examples some suitable approximations of $c_0(\theta)$ that are simpler to compute.

Suppose, for sample size $n=1$, there is a statistic $t_1\equiv T(y_1)$ that satisfies R1, i.e. it allows us to construct a valid confidence density $c_m(\theta,t_1)$. Then we can compute $c_{0}(\theta)$ based on $c_m(\theta;t_1)/L(\theta;t_1)$. First consider the case when $c_0(\theta)$ is free of the data. From the updating formula in Pawitan and Lee (2021), the 
confidence density based on the whole data is 
\begin{eqnarray}
c_{f}(\theta ;y) &\propto &c_{m}(\theta ;t_{1})L(\theta;y_1|t_1)L(\theta ;y_{2}\cdots y_{n}) \notag \\
&\propto&c_{0}(\theta )L(\theta ;t_{1})L(\theta;y_1|t_1)L(\theta ;y_{2}\cdots y_{n})  \notag \\
&=&c_{0}(\theta )L(\theta ;y_{1})L(\theta ;y_{2}\cdots y_{n})  \notag \\
&\propto &c_{0}(\theta )L(\theta ;y).  \label{eq:update} 
\end{eqnarray}
The statistic $t_{1}$ trivially exists if $y_{1}$ itself leads to a valid confidence density. Once $c_{0}(\theta )$ is available, (\ref{eq:update}) is highly convenient, since it does not require any computation of a statistic such as the MLE, its distribution or the P-value based on the whole data. More importantly, as shown in some examples below, formula (\ref{eq:update}) works even when there is no sufficient statistic from the whole data for $n>1$. This is illustrated by the general location-family model in Section~\ref{sec:location}.

When $c_{0}(\theta)$ depends on the data, it matters which $y_i$ is used to compute it. In this case the updating formula (\ref{eq:update}) is only an approximation. As long as the contribution of $\log c_0(\theta)$ to $\log c_f(\theta)$ is of order $O(1/n)$, we expect a close approximation. This is illustrated in Example~6 below.


\subsection{Simple models}

\noindent \textbf{Example 1 (continued).}
Based on $y_1$, the confidence density and the likelihood functions are
proportional: 
\begin{equation*}
c(\theta;y_1) \propto L(\theta;y_1) = 1,\ \mbox{for } \theta\in
\{y_1-1,y_1,y_1+1\}, 
\end{equation*}
so the implied prior $c_0(\theta)=1$ for all $\theta$. The full likelihood
based on $(y_1,y_2)$ is 
\begin{equation*}
L(\theta) = 1,\ \mbox{for } \theta \in \{y_{(2)}-1,y_{(1)}+1\}, 
\end{equation*}
so, the full confidence density is  $c_f(\theta) \propto L(\theta)$.
For example, if $y_1=3$ and $y_2=5$, we do have 100\% confidence that 
$\theta=4$. And if $y_1=y_2=3$, we only have 33.3\% confidence for $\theta=4$,
though we have 100\% confidence for $\theta \in \{2, 3, 4\}$.

The MLE of $\theta$ is not unique, but we can choose $\widehat{\theta}=\bar{y}$ as the
MLE. It is not sufficient, but $(\bar{y},R)$ is, so $R$ is a maximal ancillary. Indeed
the full confidence values match the conditional probabilities given the
range $R$ as previously given. Furthermore, according to Corollary 3, there is
no further relevant subset, so the confidence is epistemic. 
$\Box$\vspace{.1in}

\noindent \textbf{Example 3.} Let $y=y_1$ be a single sample
from the uniform distribution on $[\theta, \theta+1]$, where $\theta$ is a
real number. As in the previous examples, the confidence density and the
likelihood functions are 
\begin{equation*}
c_f(\theta) \propto L(\theta) = 1,\ \mbox{for } \theta\in [y-1,y]. 
\end{equation*}
For example, if $y=1.9$, then we are 100\% confident in $0.9<\theta<1.9$;
and 90\% confident in $1.0<\theta<1.9$. The coverage probability of CIs of
the form $[y-\gamma,y]$ is indeed 
\begin{equation*}
P_\theta(\theta \in \mbox{CI}) = \gamma, 
\end{equation*}
so the confidence is epistemic.

Now let's consider trying to bet the value $\lfloor \theta
\rfloor$, the largest integer smaller than $\theta$, based on observing $y$.
What price would you give to the bet that $\lfloor \theta \rfloor=\lfloor y
\rfloor = 1$? According to the confidence distribution, it should be 0.9.
But we can show that the random variable $\lfloor y \rfloor$ is Bernoulli
shifted by $\lfloor \theta \rfloor$ and with success probability $\langle
\theta \rangle \equiv \theta- \lfloor \theta \rfloor$, the fractional part
of $\theta$. For example, if $\theta=1.6$ then $\lfloor \theta \rfloor=1$
and $\langle 1.6 \rangle = 0.6$, so $\lfloor y \rfloor$ is equal to $1+0=1$
or $1+1=2$ with probabilities 0.4 and 0.6, respectively. This means that, in
general, using $\lfloor y \rfloor$ as a guess, the `coverage' probability of
being correct is 
\begin{equation*}
P_\theta\{\lfloor Y \rfloor = \lfloor \theta \rfloor\} = 1-\langle \theta
\rangle. 
\end{equation*}
This probability varies from 0 to 1 across $\theta$, not matching the
specific confidence -- such as 0.9 above -- derived from the confidence
density.

The problem is that $\lfloor y \rfloor$ is no longer a
sufficient statistic, so its marginal distribution is not fully informative.
Now, the fractional part $\langle y \rangle$ is uniform between 0 and 1 for
any $\theta$, so it is an ancillary statistic. We can show that, conditional
on $\langle y \rangle$, the distribution of $\lfloor y \rfloor$ is
degenerate: with probability 1, it is equal to $\lfloor \theta \rfloor +1$
if $\langle y \rangle < \langle \theta \rangle$; and equal to $\lfloor
\theta \rfloor$ if $\langle y \rangle > \langle \theta \rangle$. This
conditional distribution is distinct from the unconditional version, so 
$\langle y \rangle$ is relevant. Basu (1964) and Ghosh et al. (2010) used
this example as a counter-example, where conditioning by an ancillary leads
to a puzzling degenerate distribution. But actually, it is not so puzzling:
The conditional likelihood is \emph{the same} as the full likelihood, so 
$\langle y \rangle$ is maximal ancillary. This is of course as we should
expect, since $\langle y \rangle$ together with $\lfloor y \rfloor$ form the
full data $y$.

To illustrate with real numbers, for example, on observing 
$y=1.9$ we have $\langle y \rangle=0.9$, 
so $\lfloor \theta \rfloor = \lfloor y \rfloor = 1$ 
if the unknown $\langle \theta \rangle <0.9$. Now, your
betting situation is much clearer: you will bet that $\lfloor \theta \rfloor
= 1$ if you believe that $\langle \theta \rangle <0.9$, i.e. $1<\theta<1.9$.
This is exactly the same logical situation you faced before with the full likelihood and confidence density. $\Box$ \vspace{.1in}

\subsection{\label{sec:location}Location family}

Suppose $y_1,\ldots,y_n$ are an iid sample from the location
family with density 
\begin{equation*}
p_\theta(y_i) = f(y_i-\theta), 
\end{equation*}
where $f(\cdot)$ is an arbitrary but known density function, for example the
Cauchy or normal densities. Immediately, based on $y_1$ alone, the
confidence density is 
\begin{equation*}
c(\theta;y_1) = f(y_1-\theta) = L(\theta; y_1), 
\end{equation*}
so the implied prior $c_0(\theta)=1$. So, again using formula (\ref
{eq:update}), the full confidence density is 
\begin{equation}
c_f(\theta) \propto L(\theta) = \prod_{i=1}^n f(y_i-\theta).
\label{eq:location}
\end{equation}
This is a remarkably simple way to arrive at the confidence density of 
$\theta$ without having to find the MLE and its distribution.

Without further specifications, the MLE 
$T \equiv \widehat{\theta}$ is not sufficient, 
so the marginal P-value $P_\theta(T>t)$ will not
yield the full confidence. The distribution of the residuals $(y_i-\theta)$
are free of $\theta$, so the set of differences $(y_i-y_j)$'s are ancillary.
In his classic paper, Fisher (1934) showed that 
\begin{equation*}
p_\theta(\widehat{\theta}|a) = k(a)\frac{L(\theta)}{L(\widehat{\theta})}, 
\end{equation*}
where $a$ is the set of differences from the order statistics 
$y_{(1)},\ldots,y_{(n)}$. This means that the conditional likelihood based on 
$\widehat{\theta}|a$ matches the full likelihood (\ref{eq:location}), and
the confidence of CIs based on (\ref{eq:location}) will match the
conditional coverage probability. Indeed, here $(\widehat{\theta},a)$ is
sufficient and $a$ is maximal ancillary. Overall, the confidence of CIs based on (\ref{eq:location}) is epistemic. 
\vspace{.1in}

\noindent \textbf{Example 4.} Suppose that $y=(y_{1},\cdots ,y_{n})$
are i.i.d sample from the uniform distribution on $[\theta -1,\theta +1]$. Let $y_{(1)}$
and $y_{(n)}$ be ordered statistics, then $(y_{(1)},y_{(n)})$ is a sufficient
statistic. The likelihood is given by 
\begin{equation*}
L(\theta ;y)\propto I{(y_{(n)}-1\leq \theta \leq y_{(1)}+1)}
\end{equation*}
and $\mle = y_{(1)}+1$ is a MLE.
Since the uniform distribution on $[\theta -1,\theta +1]$ is location family,
we have $c_{0}(\theta )=1$ to lead the full confidence
\begin{equation*}
c_{f}(\theta)=\frac{1}{2-a}I{(y_{(n)}-1\leq \theta \leq y_{(1)}+1)}
=\frac{1}{2-a}I{(\widehat{\theta }-(2-a)\leq \theta \leq \widehat{\theta })},
\end{equation*}
where the range $a\equiv y_{(n)}-y_{(1)}$ is a maximal ancillary. $\Box$\vspace{.1in}

\subsection{Exponential family model}
Suppose the dataset $y$ is an iid sample from the exponential family model with log-density of the form 
\begin{equation}
\log p_\theta(y_i) = \sum_{j=1}^{J}h_j(\theta) t_j(y_i) - A(\theta) + c(y_i). 
\label{eq:exp-family}
\end{equation}
The MLE is sufficient if $J=1$, but not if $J>1$. In the latter case, the family is called the curved exponential family. By Theorem 1, when $J=1$ confidence statements based on the MLE will be epistemic.  (Our theory covers the continuous case in order to get exact coverage probabilities. Many important members are discrete, which is more complicated because the definition of the P-value is not unique, and the coverage probability function is guaranteed not to match any chosen confidence level. We discuss an example in Appendix A3.)

The standard evaluation of the confidence
requires the tail probability of the distribution of the MLE, which in general has no closed form formula. Barndorff-Nielsen's (1983)
approximate conditional density of the MLE $\widehat{\theta }$ is given by 
\begin{equation}
p_{\theta }(\widehat{\theta}|a)
=k |I(\widehat{\theta })|^{1/2}\frac{L(\theta )}{L(\widehat{\theta })}
+O(n^{-1}),  
\label{eq:p-formula}
\end{equation}
where the MLE is the solution of $A^{\prime }(\theta )=\sum_{i}\sum_j h'_j(\theta) t_j(y_{i})$, $a$ is the maximal ancillary and $k$ is a normalizing constant that is free of $\theta$. 
For $J=1$ and the canonical parameter $h_1(\theta)=\theta$, the ancillary is null, and the approximation leads to the right-side P-value 
\begin{equation}
P_\theta\{Z>r^{\ast }(\theta )\},\ \ r^{\ast }(\theta )\equiv r+\frac{1}{r}\log \frac{z}{r}, \label{eq:pval}
\end{equation}
where $Z$ is the standard normal variate and 
\begin{equation*}
r=\mbox{sign}(\widehat{\theta }-\theta )\sqrt{w},\ \
z=|I(\widehat{\theta })|^{1/2}(\widehat{\theta }-\theta ),
\end{equation*}
with $w=2\log \{L(\widehat{\theta })/L(\theta )\}$ and 
$I(\widehat{\theta })$ the observed Fisher information. 
\vspace{.1in}

\begin{figure}[h!]
\centerline{\includegraphics[scale=.75]{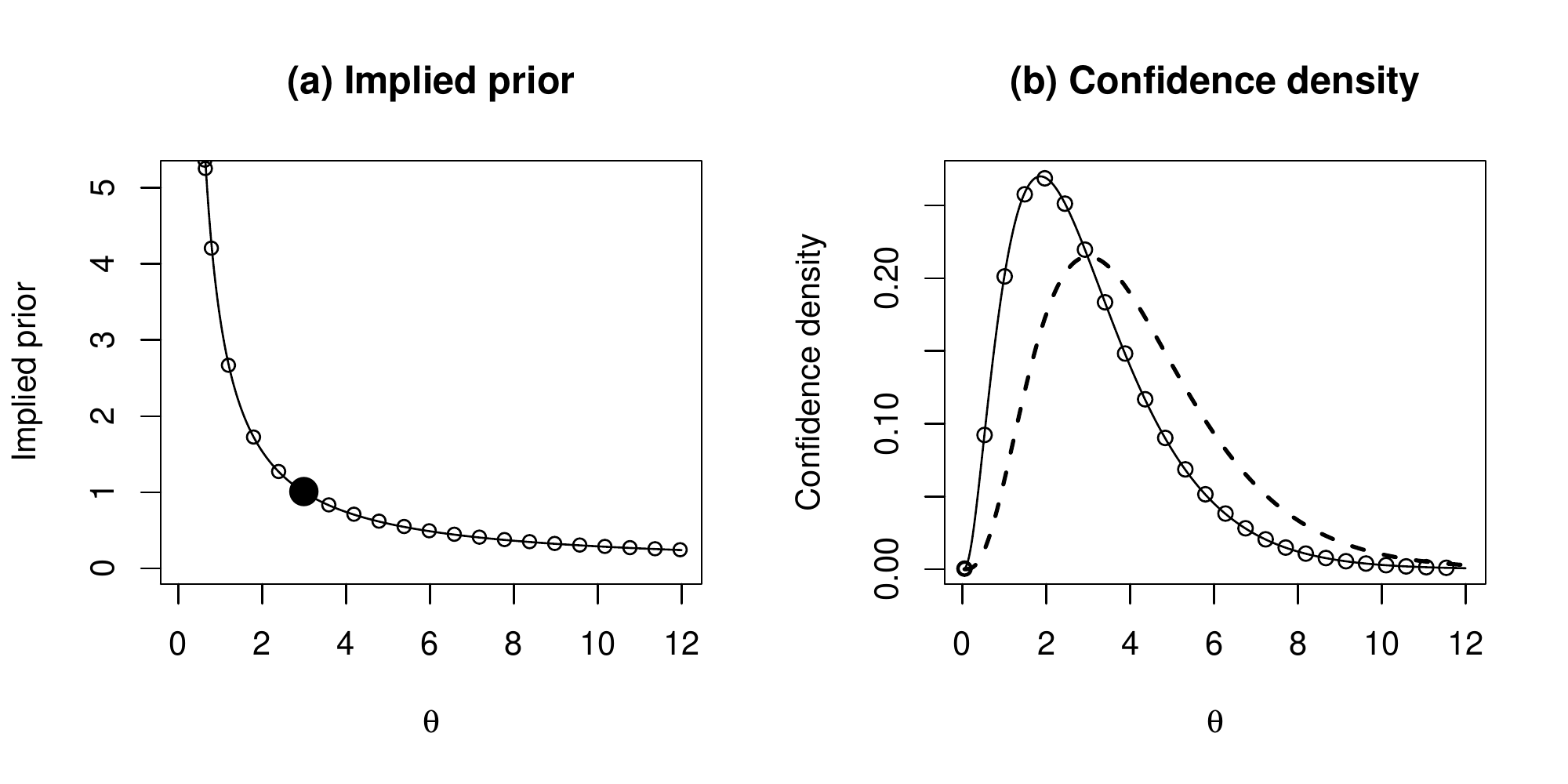}} 
\caption{\emph{(a) Implied prior of the gamma shape parameter 
$\theta $ computed using formula (\ref{eq:c0}) (solid) and from the
approximate P-value formula (\ref{eq:pval}) (circles). Both are normalized
such that they are equal to one at the MLE (black dot). (b) The confidence
densities based on a sample with size $n=5$ using formula (\ref{eq:update})
(solid) and using the approximate P-value formula (\ref{eq:pval}) (circles).
The normalized likelihood function (dashed) is also shown. }}
\label{fig1}
\end{figure}

\noindent \textbf{Example 5.}
Let $y=(y_1,\cdots,y_n)$ be an iid sample from the gamma
distribution with mean one and shape parameter $\theta$. The density is
given by 
\begin{equation*}
p_\theta(y_i) = \frac{1}{\Gamma(\theta)} \theta^\theta
y_i^{\theta-1}e^{-\theta y_i}, 
\end{equation*}
so we have an exponential family model with 
\begin{equation*}
t(y_i) = -y_i + \log y_i, \ \ A(\theta) = \log\Gamma(\theta) -
\theta\log\theta. 
\end{equation*}
To use formula (\ref{eq:update}), we first find the implied prior density using $t_1\equiv t(y_1)$ alone: 
\begin{equation*}
c_0(\theta) \propto \frac{c(\theta;t_1)}{L(\theta;t_1)},
\end{equation*}
where $c_1(\theta) = \partial \{P_\theta(T_1 \ge t_1)\}/\partial\theta$ and 
$L(\theta;t_1) = p_\theta(y_1)$. The probability $P_\theta(T_1 \ge t_1)$ is a
gamma integral, which is computed numerically. The implied prior is shown in Figure~\ref{fig1}(a). So from 
(\ref{eq:update}), we get the confidence density 
\begin{equation*}
c_f(\theta) \propto c_0(\theta)L(\theta) = c_0(\theta)\prod_{i=1}^n
p_\theta(y_i). 
\end{equation*}
For an example with $n=5$ and $\sum_i t(y_i)=-5.8791$, which corresponds to
the MLE $\widehat{\theta}=3$, the confidence density is given by the solid line
in Figure~\ref{fig1}(b). The normalized likelihood function is also shown by the dashed line, which is quite distinct from the confidence density. 

To get the marginal confidence density based on the P-value formula (\ref{eq:pval}), we need 
\begin{equation*}
w=-2n\log \Gamma (\theta )+2n\log \Gamma (\widehat{\theta })+2n(\theta \log
\theta -\widehat{\theta }\log \widehat{\theta })
+2(\theta -\widehat{\theta })\sum_{i}(\log y_{i}-y_{i})
\end{equation*}
where $\widehat{\theta }$ is the solution of 
\begin{equation*}
n\psi (\theta )-n\log \theta -n=\sum_{i}t(y_{i}),
\end{equation*}
with $\psi (\theta )\equiv \partial \log \Gamma (\theta )/\partial \theta $,
and the observed Fisher information is 
\begin{equation*}
I(\widehat{\theta })=n\{\psi ^{\prime }(\widehat{\theta })-1/\widehat{\theta 
}\}.
\end{equation*}
The circle points in Figure~\ref{fig1}(b) are the marginal confidence density based on the same sample above. As expected, this tracks almost exactly the one given by formula (\ref{eq:update}). The corresponding implied prior based on $c_{m}(\theta)/L(\theta)$ is given in Figure~\ref{fig1}(a), also closely matching the implied prior based on $n=1$. $\Box$ \vspace{.1in}

\begin{figure}[h!]
{\normalsize \centerline{\includegraphics[scale=.651]{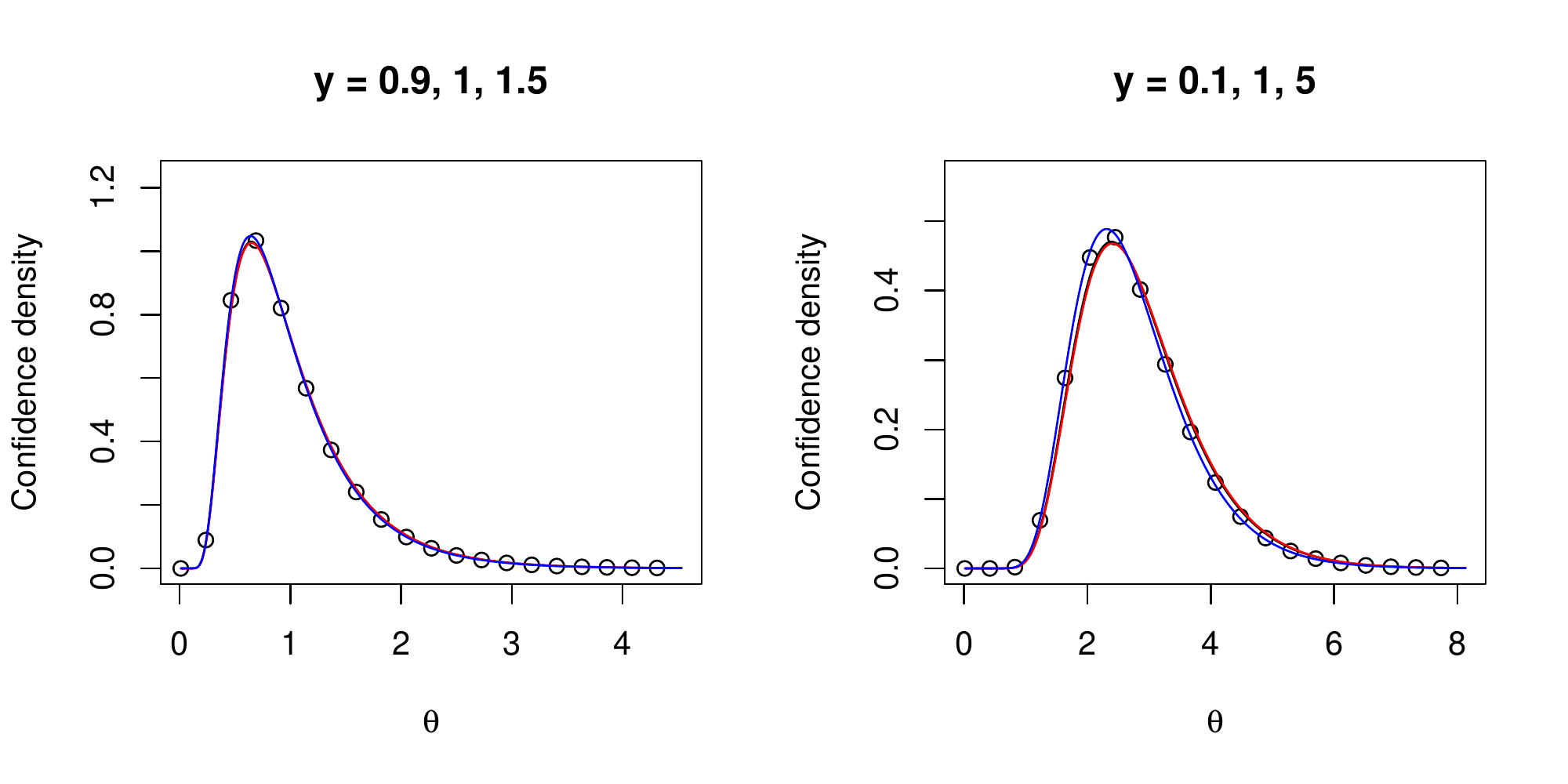}}
} \caption{\label{fig:normal-mumu}\emph{Confidence densities $c_{fi}(\theta)$ (solid) and $c_m(\theta)$ (circles) for the $N(\theta,\theta)$ model. The former is actually drawn three times, corresponding to three different versions of the implied prior based on $y_i$ for $i=1, 2, 3$.  (a) Based on $y=(0.9, 1, 1.5)$ (b) Based on $y= (0.1, 1, 5).$ }}
\end{figure}

\noindent \textbf{Example 6.}
This is an example where $c_0(\theta)$ is data dependent. Let $y=(y_1,\ldots,y_n)$ be iid sample from $N(\theta,\theta)$ for
$\theta>0$. The log-density is given by
$$
\log p_\theta(y) = -\frac{n}{2}\log(2\pi \theta) -\half (\sum_i y_i^2/\theta - 2\sum_i y_i + n \theta),
$$
so this is a regular exponential family with sufficient statistic $T(y) = \sum_i y_i^2$. The marginal confidence density $c_m(\theta)$ can be computed based on the non-central $\chi^2$ distribution for $T(y)$. For $n=1$,  $T(y_1) = y_1^2$ is sufficient, and
\begin{eqnarray*}
C(\theta;y_1) &=& P_\theta(Y_1^2>y_1^2)\\
     &=&   1-\Phi\left(\frac{|y_1|-\theta}{\sqrt{\theta}}\right) + \Phi\left(\frac{-|y_1|-\theta}{\sqrt{\theta}}\right)\\
c(\theta;y_1) &=&\half \phi\left(\frac{|y_1|-\theta}{\sqrt{\theta}}\right)
\left(\frac{|y_1|}{\theta\sqrt{\theta}}+\frac{1}{\sqrt{\theta}}\right)
 + \half \phi\left(\frac{-|y_1|-\theta}{\sqrt{\theta}}\right)
\left(\frac{|y_1|}{\theta\sqrt{\theta}}-\frac{1}{\sqrt{\theta}}\right)\\
L(\theta; y_1) &=& \phi\left(\frac{y_1-\theta}{\sqrt{\theta}}\right),
\end{eqnarray*}
so the implied prior is data-dependent. This means that the full confidence density depends on which $y_i$ is used to compute the implied prior:
$$
c_{fi}(\theta) = c_0(\theta;y_i) L (\theta;y)
$$
In Figure~\ref{fig:normal-mumu}, for $n=3$, we compare $c_{fi}(\theta)$ using three different versions of $c_0(\theta)$ based $y_i$ for $i=1, 2, 3$. These are also compared with the marginal confidence $c_m(\theta)$. As shown in the figure, even for such a small dataset, the effects of the data dependence in this case are negligible. $\Box$ \vspace{.1in}

\noindent \textbf{Example 7.} 
This example from a curved exponential family is used to illustrate complex models, 
where standard results sometimes fail. Let $y_1,\ldots,y_n$ be iid sample from $N(\theta,\theta^2)$ for $\theta > 0$. First consider the confidence distribution based on $y_1$, 
$$
C_m(\theta;y_1) = P_\theta(Y_1 \geq y_1) = 1- \Phi\left(\frac{y_1-\theta}{\theta}\right).
$$
We can see immediately that if we use $y_1$ as the statistic, 
the term inside the bracket converges to $-1$ as $\theta\rightarrow \infty$, 
so the confidence distribution goes to $1-\Phi(-1) = 0.84$. 
Hence $y_1$ does not satisfy the regularity condition R1.

Here $(\sum y_i^2, \sum y_i)$ is minimal sufficient, and the likelihood function is
$$
L(\theta; y) = {(2\pi\theta^2)^{-n/2}} \exp \left( -\sum y_i^2 / 2\theta^2 + \sum y_i /\theta - n/2 \right).
$$
The MLE is given by
$$
\mle = \mle (y) = \frac{-\sum y_i + \sqrt{(\sum y_i)^2 + 4 \sum y_i^2}}{2n}
$$
with a maximal ancillary
\begin{align*}
A(y) = \frac{\sum y_i}{\sqrt{\sum y_i^2}}.
\end{align*}
In terms of $(\mle,a)\equiv (\mle(y),A(y)) $, the likelihood is
\begin{align*}
L(\theta; y) = {(2\pi\theta^2)^{-n/2}} 
\exp \left( -\frac{(b + 2n)\mle^2}{4\theta^2} 
+ \frac{b\mle}{2\theta} - \frac{n}{2} \right),
\end{align*}
where $b\equiv a^2+a\sqrt{a^2+4n}$.

Now,  denote the MLE and the ancillary based only on $y_1$ by $\mle_1=\that(y_1)$ and $a_1$. The conditional confidence distribution based on $\that_1|a_1$ is
\begin{align*}
C_c(\theta; \mle_1|a_1) = P_\theta (\mle(Y_1) \geq \mle_1 | a_1 )
= \frac{1}{\Phi(a_1)} \Phi \left( \frac{-a_1-\sqrt{5}}{2} \frac{\mle_1}{\theta} + a_1 \right),
\end{align*}
which is now a valid confidence distribution, with density
\begin{align*}
c_c(\theta; \mle_1|a_1) = \frac{\partial}{\partial \theta} C_c(\theta;\mle_1 | a_1)
= \frac{1}{\Phi(a_1)} \frac{a_1+\sqrt{5}}{2} \frac{\mle_1}{\theta^2} \phi \left( \frac{-a_1-\sqrt{5}}{2} \frac{\mle_1}{\theta} + a_1 \right).
\end{align*}
The implied prior is $c_0(\theta; \mle_1 |a_1) \propto c_c(\theta; \mle_1|a_1)/L(\theta; y_1) \propto \theta^{-1}$. The updating rule gives the full confidence density
\begin{align} \label{eq:mu2_cf}
c_f(\theta; y) \propto c_0(\theta) L(\theta; y)
&\propto \frac{1}{\theta^{n+1}}\exp \left( \sum \frac{y_i^2}{2\theta^2} + \sum \frac{y_i}{\theta} \right) \nonumber \\
&\propto \frac{1}{\theta^{n+1}} \exp \left( -\frac{(b + 2n)\mle^2}{4\theta^2} 
+ \frac{b\mle}{2\theta} \right).
\end{align}
In fact, in Appendix A2 we show that, even though it is not sufficient, $\that_1$ still leads to a valid confidence distribution, and the implied prior based on $c_m(\theta; \mle_1)/L(\theta; \mle_1)$ is the same $c_0(\theta) = 1/\theta$.

For completeness, Appendix A2 also shows the conditional confidence density derived using Barndorff-Nielsen's formula (\ref{eq:p-formula}), showing that we end up with the same implied prior. Instead here we use the exact result from Hinkley (1977). He derived the exact conditional density of $w=\theta^{-1}\sqrt{\sum y_i^2}$,
$$
p(w|a) = w^{n-1} \exp \{ -(w - a)^2 /2 \} / I_{n-1} (a)
$$
where $I_{n-1}(a)=\int_0^{\infty} x^{n-1} \exp\{-(x-a)^2/2\}dx$. Let $T(y)=\sqrt{\sum y_i^2}$, then we have
$$
C_c(\theta; t|a) = P_\theta (T\geq t |a) = P(W\geq w |a) = 1 - F_a(w)
$$
where $F_a(w) = \int p(w|a) dw$. Then the confidence density becomes
$$
c_c(\theta; t|a) = - \frac{\partial F_a(w)}{\partial\theta} = p(w=t/\theta|a) \ t/\theta^2 = p_\theta(t|a) \ t/\theta
$$
so that the implied prior becomes $c_0(\theta; t|a) \propto c_c(\theta; t|a)/L(\theta; t|a) \propto 1/\theta$,
which is the same with the results from $\mle|a$. Thus, we have
$$
c_c(\theta; t|a) = c_c(\theta; \mle|a) = c_f(\theta;y).
$$

In Appendix 2, it is also shown that $c_m(\theta; \mle_1,\cdots,\mle_n)$ with $\mle_i = \mle(y_i)$ is a valid confidence density because $C_m(\theta\in \CI(\mle_1,\cdots,\mle_n)) = P_\theta(\theta \in \CI(\mle_1,\cdots,\mle_n))$. However, it is not epistemic because it does not use the full likelihood, so there is a loss of information.  

\begin{figure}[h!]
{\normalsize \centerline{\includegraphics[scale=.8]{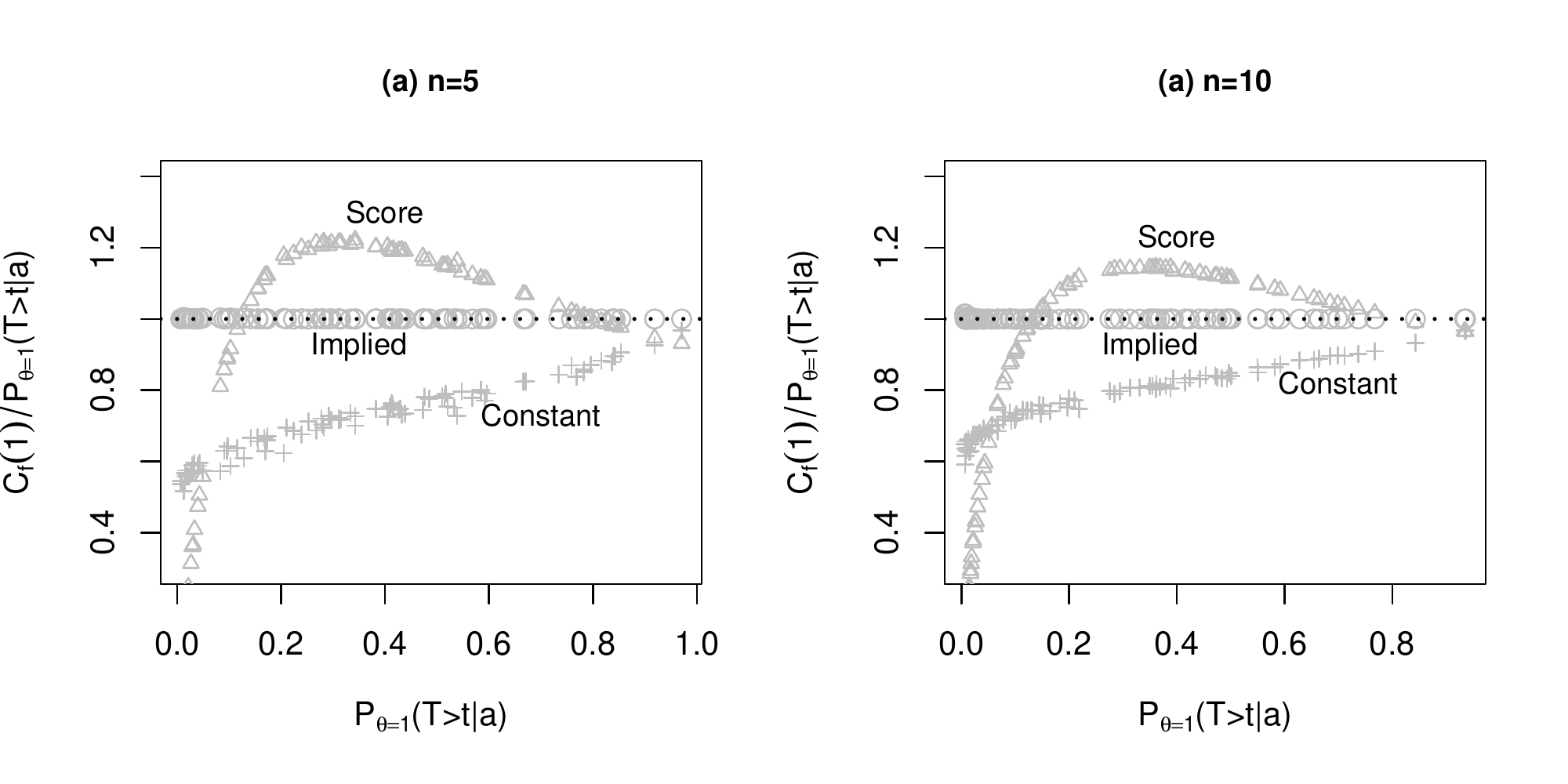}}
}
\caption{\label{fig:curved}\emph{Example from the $N(\theta,\theta^2)$ model. 
In each panel, the x-axis is the exact conditional P-value $P_\theta(T>t|a)$ given the ancillary $a$ for testing $H_0$: $\theta=1$. The y-axis is the full confidence $C_f(1) \equiv \int_{\theta<1} c_0(\theta)L(\theta;y)/m(y) d\theta$, 
using the constant prior $c_0(\theta)= 1$ ('+' symbols) and the implied prior $c_0(\theta)= 1/\theta$ (circles).
Also shown is the corresponding P-value from the score test using Fisher's observed information (triangles). 
(a) For $n=5$ and (b) for $n=10$. To show the quality of the approximation for small P-values, the y-axis is expressed as a ratio.}}
\end{figure}

As numerical illustrations, we compare the exact conditional P-value $P_\theta(T>t|a)$ for testing $H_0$: $\theta=1$, the corresponding full confidence $C_f(\theta)$ at $\theta=1$ and the P-value based on the score test. The latter was computed using the observed Fisher information,
suggested by Hinkley (1977) as having good conditional properties. In Figure~\ref{fig:curved}(a), we generate 100 datasets with $n=5$ from $N(\theta,\theta^2)$ at $\theta=1.2$.
The full confidence $C_f(\theta<1)$ is computed using the implied prior $c_0(\theta) \propto 1/\theta$, and a constant prior $c_0(\theta) \propto 1$. Panel (b) shows the result for $n=10$.
The full confidence with the implied prior $c_0(\theta)\propto 1/\theta$ agrees with the exact conditional probability.
The use of a non-implied prior $c_0(\theta) \propto 1$ cannot match the exact conditional probability. While the score test maintains the average conditional probability, it has a poor conditional property in these small samples. $\Box$ \vspace{.1in}

\section{\label{sec:disc}Discussion}

We have described the Dutch Book argument to establish the epistemic confidence that is meaningful for an observed confidence interval. Fisher tried to achieve the same purpose with the fiducial probability, but the use of the word 'probability' had generated much confusion and controversies, so the concept of fiducial probability has been practically abandoned. However,
the confidence concept is mainstream, although it comes
with a frequentist interpretation only, so it applies not  to the observed interval but to the procedure.
The confidence may not be a probability but an extended likelihood (Pawitan and Lee, 2021),
whose ratio is meaningful in hypothesis testing and statistical inferences (Lee and Bj{\o}rnstad, 2013).
The extended likelihood is logically distinct from the classical likelihood. It is well known that we cannot use the likelihood directly for inference, except in special circumstances such as the normal case. Our results show that we can turn a classical likelihood into a confidence density by multiplying it with an implied prior. Furthermore, we get epistemic confidence by establishing the absence of relevance subsets.

It might appear that in trying to get epistemic confidence using the Dutch Book argument, we are simply recreating the Bayesian framework. But this is not the case. Bayesian subjective probability is established by an internal consistency requirement of a coherent betting strategy. If you're being
inconsistent, then an external agent can construct a Dutch Book against you. To avoid the Dutch Book, you have to use a (additive) probability measure. Crucially, in this Bayesian version, there is no reference to the betting market, which according to the fundamental theorem of asset pricing (Ross, 1976) will settle prices based on objective probabilities. So, in our setting, it is still possible to construct a market-linked Dutch Book against an internally consistent subjective Bayesian who ignores the objective probability.

Extra principles in the subjective probability framework have been proposed to deal the mismatch between the subjective and objective
probabilities. For example, in Lewis's (1980) `Principal Principle'
\begin{equation*}
P_s\{A| \mbox{Chance}(A)=x\} = x,
\end{equation*}
where $P_s$ denotes the subjective probability and `Chance' the objective probability.
So, the Principle simply declares that the subjective probability must be set to be
equal to the objective probability, if the latter exists. Our Dutch Book argument can be used to justify the Principle, so the principle does not have to come out of the blue with no logical motivation.
However, it is worth noting that epistemic confidence is \emph{not simply set equal to probability}. Instead, it is the consequence of
a theorem that establishes no relevant subset in order to avoid the
Dutch Book. In our setup, the frequentist probability applies to a market involving a large number of independent players. Moreover, the rational personal betting price is no longer `subjective', for example, in the choice of the prior. Thus, the conceptual separation of the personal and the market prices allow both epistemic and frequentist meaning of confidence.

Our use of money and bets to define epistemic confidence has some echoes in Shafer and Vovk (2001)'s game-theoretic foundation of probability, an ambitious rebuilding of probability without measure theory. However, their key concept is a sequential game between two players. The word `sequential' clearly implies that the game is not meant to involve a risk-free profit from a single transaction that we want in a Dutch Book. Our usage of probability is fully within the Kolmogorov axiomatic system, and we make a clear distinction between probability and confidence.

Conditional inference (Reid 1995) has traditionally been the main area of statistics that tries to address the epistemic content of
confidence intervals. The theory on ancillary statistics and relevant subsets has grown as a result (Basu, 1955, 1954; Ghosh, 2010).
Conditioning on ancillary statistics is meant to make the inference `closer' to the data at hand, but the proponents of conditional
inference only go half-way to the end goal of epistemic confidence that Fisher wanted. The general lack of unique maximal ancillary is a great stumbling block, where it is then possible to come up with distinct relevant subsets with distinct conditional coverage probabilities. This raises an unanswerable question: What is then the `proper' confidence for the observed interval? Our logical tool of the betting market overcomes this problem -- in this case, the market cannot settle in an unambiguous price. But Corollary 3 still hold in the sense that as an individual, you're still protected from the Dutch Book. We discuss this further with an example in Appendix 4.

Schweder and Hjort (2016) and Schweder (2018)
have been strong proponents of interpreting confidence as `epistemic probability.' We are in general agreement with their sentiment, but it is unclear which version of probability this is. The only established and accepted epistemic probability is the Bayesian probability, but in their writing, the confidence concept is clearly non-Bayesian. Our use of the Dutch Book defines normatively the epistemic value of the confidence while staying within the non-Bayesian framework.

In logic and philosophy, the relevant subset problem is known as the `reference class problem' (Reichenbach, 1949; Hajek, 2007). Venn (1876) in his classic book already recognized the problem: `It
is obvious that every individual thing or event has an indefinite
number of properties or attributes observable in it, and might
therefore be considered as belonging to an indefinite number of
different classes of things...', and this affects the probability assignment. We solve this problem by setting a limit to the amount of information as given by the observed data $y$ generated by the probability model $p_\theta(y)$. This approach is in line with the theory of market equilibrium, where it is assumed that information is limited and available to all players. It is not possible to have an equilibrium -- hence accepted prices -- if information is indefinite, or when the players know that they have access to different pieces of information.


We have limited our current paper to the one parameter case. Following the proof of Theorem 1, the same result actually holds in the multi-parameter case as long as we consider bounded confidence regions satisfying Condition R4. The problem arises for a marginal parameter of interest, which might implicitly assume unbounded regions for the nuisance parameters. For example, in the normal model, with both mean $\mu $ and variance $\sigma^{2}$ unknown, the standard CI for the mean implicitly employs an unbounded
interval for the variance:
\begin{equation*}
\mu \in \left( \bar{y}-t_{\alpha /2}\frac{s}{\sqrt{n}},\bar{y}+t_{\alpha /2}
\frac{s}{\sqrt{n}}\right) \ \text{and}\ \sigma ^{2}\in (0,\infty )
\end{equation*}
where $t_{\alpha /2}$ is the quantile of $t$ distribution with $(n-1)$
degrees of freedom. Thus, the theorem cannot guarantee the epistemic property of the $t$-interval. In fact, Buehler and Feddersen (1963) and Brown (1967) showed the existence of relevant subsets for the $t$-interval. Extensions of our epistemic confidence results to this problem are of great interest.

\section*{References}

\begin{description}
\item \textsc{Arrow, K. J. and Debreu, G.} (1954). Existence
of an equilibrium for a competitive economy. \emph{Econometrica} 
\textbf{22}, 265–-290.

\item \textsc{Barndorff-Nielsen O.} (1983). 
On a formula for the distribution of the maximum likelihood estimator. \emph{Biometrika}, \textbf{70}: 343-365


\item \textsc{Berger J.O. and Wolpert R.L.} (1988). \emph{The likelihood principle.} 2nd Edition. Institute of Mathematical Statistics,
Hayward, CA 

\item \textsc{Birnbaum A.}\ (1962). On the foundation of
statistical inference. \emph{Journal of the American Statistical Association,}
\textbf{57}, 269--326.

\item \textsc{Bj{\o }rnstad J.F.}\ (1996). On the
generalization of the likelihood function and likelihood principle. \emph{%
Journal of the American Statistical Association,} \textbf{91}, 791--806. 

\item \textsc{Brown, L.}\ (1963). The conditional level of
Student's t test. \emph{The Annals of Mathematical Statistics,} \textbf{38},
1068--1071.

\item \textsc{Buehler R.J.} (1959). Some Validity Criteria for Statistical Inferences. 
\emph{The Annals of Mathematical Statistics,}
\textbf{30}, 845--863

\item \textsc{Buehler R.J. and Feddersen A.P.}\ (1963). Note
on a Conditional Property of Student's $t$. \emph{The Annals of Mathematical Statistics,} \textbf{34}, 1098--1100.



\item \textsc{de Finetti, B.} (1931). On the subjective meaning of probability. English translation in de
Finetti (1993), \emph{Induction and Probability,} pp.\ 291–321. Clueb: Bologna. 


\item \textsc{Fisher R.A.}\ (1930). Inverse probability. \emph{%
Proceedings of the Cambridge Philosophical Society,} \textbf{26}, 528--535.

\item Fisher R.A.\ (1933). The concepts of inverse probability and fiducial probability referring to unknown parameters. 
\emph{Proceedings of the Royal Society of London,} \textbf{139}, 343--348.

\item \textsc{Fisher R.A.}\ (1934). Two new properties of
mathematical likelihood. \emph{Proceedings of the Royal Society of London}, 144A, 285.

\item \textsc{Fisher R.A.}\ (1950). 
\emph{Contributions to mathematical statistics} 
(pp. 35.173a). Wiley.

\item \textsc{Fisher R.A.}\ (1958). The nature of probability. 
\emph{Centennial Review,} \textbf{2}, 261--274.

\item \textsc{Fisher R.A.}\ (1973). \emph{Statistical methods and scientific inference}, 3rd Edition. New York: Hafner. 

\item \textsc{Ghosh M., Reid N.\ and Fraser D.A.S.} (2010).
Ancillary statistics: a review. \emph{Statistica Sinica} 20, 1309--1332.

\item \textsc{Hajek A}. (2007). The reference class problem is
your problem too. \emph{Synthese}, \textbf{156}: 563--585.

\item \textsc{Lancaster H.O.} (1961). Significance tests in
discrete distributions. \emph{Journal of the American Statistical Association,} 
\textbf{56}: 223--234.

\item \textsc{Lee Y. and Bj{\o}rnstad J.F.} (2013). 
Extended likelihood approach to large-scale multiple testing.
\emph{Journal of the Royal Statistical Society. Series B (Statistical Methodology)}, \textbf{75}, 553-575.

\item \textsc{Lee Y., Nelder J.A. and Pawitan Y. } (2017). 
\emph{Generalized linear models with random effects: unified analysis via
H-likelihood} (2nd ed.). CRC Press.

\item \textsc{Lewis, D.} (1980). A subjectivist’s guide to objective chance. In \emph{Ifs} (pp. 267-297). Springer, Dordrecht.

\item \textsc{Lindley D.V.} (1958). Fiducial distributions
and Bayes' theorem. \emph{Journal of the Royal Statistical Society. Series B
 (Statistical Methodology)}, \textbf{20}, 102--107.

\item \textsc{Pawitan Y.} (2001). \emph{In all likelihood:
Statistical modelling and inference using likelihood}. Oxford University
Press, Oxford, UK.

\item \textsc{Pawitan Y. and Lee Y.} (2021). Confidence as
likelihood. \emph{Statistical Science.} To appear.

\item \textsc{Ramsey, F.} (1926). \emph{Truth and probability}. In \emph{Foundations of Mathematics and other Logical Essays}. London: K.\ Paul, Trench, Trubner and Co. Reprinted in H.E. Kyburg and H.E. Smokler (eds.) (1980), \emph{Studies in Subjective Probability}, 2nd edn (pp. 25-–52) New York: Robert Krieger.

\item \textsc{Reichenbach H.} (1949). \emph{The Theory of
Probability}. University of California Press. 

\item \textsc{Reid N.} (1995). The roles of conditioning in
inference. \emph{Statistical Science,} \textbf{10}, 138--157.

\item \textsc{Robinson G. K.}. (1979). Conditional properties
of statistical procedures. \emph{The Annals of Statistics}, \textbf{7},
742--755.

\item \textsc{Ross, S.} (1976), Return, risk and arbitrage.
In: I. Friend and J. Bicksler (eds.), Studies in Risk and Return. Cambridge,
MA: Ballinger


\item \textsc{Schweder T. and Hjort N.L.} (2016). \emph{%
Confidence, Likelihood, Probability.} Cambridge University Press, Cambridge,
UK.

\item \textsc{Schweder T.} (2018). Confidence is epistemic
probability for empirical science. \emph{Journal of Statistical Planning and
Inference} \textbf{195}, 116--125.

\item \textsc{Venn, J}. (1876). \emph{The Logic of Chance}
(2nd ed.). Macmillan and Co.


\end{description}
 
\section*{Appendix}

\subsection*{A1. Proof of Theorem 1}
One main consequence of condition R1 is that we start with a proper confidence density. However, as stated after the statement of the theorem, instead of starting with $T$, for more generality, our proof would allow an arbitrary function $c_0(\theta,y)$ that satisfies R2 and R3, as long as the resulting full confidence $c(\theta;y)\propto c_0(\theta,y)L(\theta;y)$ is proper.

Throughout this part, given a fixed $\gamma \in (0,1)$, let $\CI_\gamma(y)$ denote the $\gamma$-level interval 
such that
$$
C(\theta\in \mbox{CI}_\gamma(y))=\int_{CI_\gamma(y)} c(\theta;y) d\theta = \gamma
\quad \text{for any } y \in \mathcal{Y},
$$
where $\mathcal{Y}$ is the sample space of $Y$ and let $R$ denote a relevant subset such that
\begin{equation}
P_\theta (\theta \in CI_\gamma(Y)|Y\in R) \geq \gamma+\epsilon \quad \text{for any }
\theta \in \Theta  \label{eq:relevant}
\end{equation}
or negatively, 
\begin{equation*}
P_\theta (\theta \in CI_\gamma(Y)|Y\in R) \leq \gamma-\epsilon \quad \text{for any }
\theta \in \Theta
\end{equation*}
for some $\epsilon>0$ where $\Theta$ is the parameter space of $\theta$.\vspace{.1in}

\noindent \textbf{Lemma.} Suppose that $R$ is a relevant subset for an interval $\CI_\gamma(Y)$. Let $\{ B_k \}$ be a countable partition of $R$,
so that $R= \bigsqcup_{k=1}^{\infty} B_k$, then there exists a relevant
subset $B_i \in \{ B_k \}$.

\noindent \textbf{Proof.} First consider the positively biased case. If $%
P_\theta (\theta \in CI_\gamma(Y) | Y\in B_k) < \gamma + \epsilon$ for all $k$,
then we have 
\begin{align*}
P_\theta (\theta \in CI_\gamma(Y)|Y\in R) &= \frac{P_\theta (\theta \in CI_\gamma(Y) \ \&
\ Y\in R)}{P_\theta (Y\in R)} \\
&= \frac{\sum_{k=1}^{\infty} P_\theta(\theta \in CI_\gamma(Y) \ \& \ Y\in B_k)}{%
\sum_{k=1}^{\infty} P_\theta (Y\in B_k)} \\
&< \frac{\sum_{k=1}^{\infty} (\gamma+\epsilon) P_\theta (Y\in B_k)}{%
\sum_{k=1}^{\infty} P_\theta (Y\in B_k)} = \gamma + \epsilon,
\end{align*}
which is a contradiction to \eqref{eq:relevant}. Therefore, there should be $%
B_i \in \{ B_k \}$ such that 
\begin{align*}
P_\theta (\theta \in CI_\gamma(Y) | Y\in B_i) \geq \gamma + \epsilon
\end{align*}
is a relevant subset. Similarly, supposing $P_\theta (\theta \in CI_\gamma(Y) |
Y\in B_i) \leq \gamma - \epsilon$ leads to the corresponding proof for a
negatively biased relevant subset. $\Box$ \vspace{.1in}

\noindent \textbf{Proof of Theorem 1.} By definition of $c(\theta;y)$ above, there is a non-negative function $m(y)$ such that  $c_0(\theta;y) = m(y) c(\theta;y)/L(\theta;y)$.
First we show that the existence of positively biased relevant subset $R$
leads to a contradiction. For an arbitrary $\gamma \in (0,1)$, 
suppose that there exists $\epsilon>0$ such that 
\begin{equation*}
P_\theta (\theta \in \mbox{CI}_\gamma(Y)|Y\in R) \geq \gamma+\epsilon \quad 
\text{for any } \theta \in \Theta.
\end{equation*}
Let $\delta$ and $\xi$ be arbitrary numbers satisfying $\gamma
(1-\gamma-\epsilon) / (\gamma+\epsilon) < \delta < (1-\gamma)$ and 
\begin{align} \label{eq:xi}
0< \xi < \frac{1}{2} \log \left( \frac{(\gamma+\epsilon) (\gamma+\delta)}{%
\gamma} \right).
\end{align}
By the uniform continuity of $\log c_0(\theta; y)$, there exists $d>0$ such
that 
\begin{align*}
&\| y_1 - y_2 \| < d \Rightarrow |\log c_0(\theta;y_1) - \log
c_0(\theta;y_2)| < \xi
\end{align*}
Let $\{ A_k \}$ be a partition of $\mathbb{R}^n$ divided by $n$-dimensional
grid with spacing $\sqrt[n]{d}$ where $n$ is the sample size. Define $B_k =
A_k \cap R$, then $\{ B_k \}$ becomes a countable partition of relevant
subset $R$. By the Lemma, there exists a relevant subset $B_i \in \{ B_k \}$
such that 
\begin{equation*}
P_\theta (\theta \in \mbox{CI}_\gamma(Y)|Y\in B_i) \geq \gamma+\epsilon
\quad \text{for any } \theta \in \Theta.
\end{equation*}
It can be written by the integration form, 
\begin{equation} \label{eq:int1}
\int_{B_i} I_{(\theta \in CI_\gamma(y))} f_\theta (y) dy \geq \int_{B_i}
(\gamma+\epsilon) f_\theta (y) dy \quad \text{for any } \theta \in \Theta.
\end{equation}
Since $\| y_1 - y_2 \|<d$ for any $y_1, y_2 \in B_i$, from the regularity
condition R3 we have 
\begin{align} \label{ineq:xi}
e^{-\xi} < \frac{c_0(\theta;y_1)}{c_0(\theta;y_2)} < e^{\xi } \quad \text{%
for any } y_1, y_2 \in B_i.
\end{align}
Let $\mbox{CI}_{\gamma+\delta} (Y) \supset \mbox{CI}_\gamma(Y)$ be the $%
(\gamma+\delta)$-level confidence interval. By the regularity condition R4,
there exists a compact set $\Theta_i \subseteq \Theta$ such that $\mbox{CI}%
_{\gamma+\delta} (y) \subseteq \Theta_i$ for any $y \in B_i$. Then, 
\begin{align} \label{ineq:delta}
\int_{\Theta_i} c(\theta;y)d\theta \geq \int_{CI_{\gamma+\delta}(y)}
c(\theta;y) d\theta = \gamma +\delta \quad \text{for any } y\in B_i.
\end{align}
Now let $y^*$ be an arbitrary observation in $B_i$ and integrate %
\eqref{eq:int1} over $\Theta_i$ with the $c_0(\theta; y^*)$, then we have 
\begin{equation} \label{eq:int2}
\int_{\Theta_i} \int_{B_i} I_{(\theta \in CI_\gamma(y))} f_\theta(y) dy
c_0(\theta;y^*) d\theta \geq \int_{\Theta_i} \int_{B_i} (\gamma+\epsilon)
f_\theta (y) dy c_0(\theta;y^*) d\theta.
\end{equation}
By the regularity condition R2, 
\begin{align*}
0 & \leq \int_{\Theta_i} \int_{B_i} I_{(\theta \in CI_\gamma (y))} f_\theta
(y) dy c_0(\theta;y^*) d\theta \leq \int_{\Theta_i} \int_{B_i} f_\theta (y)
dy c_0(\theta;y^*) d\theta \\
&\leq \int_{\Theta_i} P_\theta (Y \in B_i) c_0(\theta; y^*) d\theta \leq
\int_{\Theta_i} c_0(\theta; y^*) d\theta < \infty.
\end{align*}
Thus, the order of integration in \eqref{eq:int2} can be exchanged by
Fubini's theorem to lead 
\begin{align} \label{eq:exchange}
\int_{B_i} \int_{\Theta_i} I_{(\theta \in CI(y))} f_\theta(y) c_0(\theta;
y^*) d\theta dy \geq \int_{B_i} \int_{\Theta_i} (\gamma+\epsilon) f_\theta
(y) c_0(\theta; y^*) d\theta dy.
\end{align}
Note here that the inequality \eqref{ineq:xi} implies that 
\begin{align*}
e^{-\xi}c(\theta;y)m(y) < f_\theta(y) c_0(\theta;y^*) = f_\theta(y)
c_0(\theta;y ) \frac{c_0(\theta;y^*)}{c_0(\theta;y)} <
e^{\xi}c(\theta;y)m(y).
\end{align*}
Then the left-hand-side of \eqref{eq:exchange} becomes 
\begin{align*}
\int_{B_i} \int_{\Theta_i} I_{(\theta \in CI_\gamma(y))} f_\theta(y)
c_0(\theta; y^*) d\theta dy & = \int_{B_i} \int_{CI_\gamma (y)} f_\theta(y)
c_0(\theta; y^*) d\theta dy \\
&< \int_{B_i} e^{\xi} \int_{CI_\gamma(y)} c(\theta; y) d\theta m(y)dy =
e^{\xi} \gamma \int_{B_i} m(y) dy
\end{align*}
and the right-hand-side of \eqref{eq:exchange} becomes 
\begin{align*}
\int_{B_i} \int_{\Theta_i} (\gamma+\epsilon) f_\theta(y) c_0(\theta;y^*)
d\theta dy & >\int_{B_i} (\gamma+\epsilon) e^{-\xi} \int_{\Theta_i}
c(\theta;y) d\theta m(y) dy \\
& \geq e^{-\xi} (\gamma+\epsilon) (\gamma+\delta) \int_{B_i} m(y) dy.
\end{align*}
by the inequality \eqref{ineq:delta}. Thus we have 
\begin{equation*} \label{eq:final}
e^{\xi} \gamma > e^{-\xi} (\gamma+\epsilon) (\gamma+\delta),
\end{equation*}
which is a contradiction to \eqref{eq:xi}.

For negatively biased cases, define $\xi$ by an arbitrary number 
\begin{align} \label{eq:xi_neg}
0 < \xi < \frac{1}{2} \log \left( \frac{\gamma}{\gamma-\epsilon} \right)
\end{align}
and replace the inequality \eqref{eq:exchange} by 
\begin{align*}
\int_{B_i} \int_{\Theta_i} I_{(\theta \in CI(y))} f_\theta(y) c_0(\theta;
y^*) d\theta dy \leq \int_{B_i} \int_{\Theta_i} (\gamma-\epsilon) f_\theta
(y) c_0(\theta; y^*) d\theta dy.
\end{align*}
Then the left-hand-side becomes 
\begin{align*}
\int_{B_i} \int_{\Theta_i} I_{(\theta \in CI(y))} f_\theta(y) c_0(\theta;
y^*) d\theta dy > e^{-\xi}\gamma \int_{B_i} m(y)dy
\end{align*}
and the right-hand-side becomes 
\begin{align*}
\int_{B_i} \int_{\Theta_i} (\gamma-\epsilon) f_\theta (y) c_0(\theta; y^*)
d\theta dy <e^{\xi} (\gamma - \epsilon) \int_{B_i} m(y)dy
\end{align*}
which lead to a contradiction to \eqref{eq:xi_neg}. $\Box$

\subsection*{A2. Curved exponential family: $N(\theta,\theta^2)$}

We give more details of the $N(\theta,\theta^2)$ model. Denote the MLE based only on $y_1$ by $\mle(y_1)=\mle_1$. The confidence distribution is
\begin{align*}
C_m(\theta; \mle_1) 
&= P_\theta (\mle(Y_1) \geq \mle_1) 
= P_\theta \left( \frac{-Y_1 + \sqrt{5 Y_1^2}}{2}\geq \mle_1 \right) \\
&= P_\theta \left( Y_1 \geq \frac{2\mle_1}{\sqrt{5}-1} \right) +
P_\theta \left( Y_1 \leq \frac{-2\mle_1}{\sqrt{5}+1} \right) \\
&= 1 - \Phi \left( \frac{1+\sqrt{5}}{2} \frac{\mle_1}{\theta} -1 \right)
+ \Phi \left( \frac{1-\sqrt{5}}{2} \frac{\mle_1}{\theta} -1 \right).
\end{align*}
Then $\lim_{\theta \rightarrow \infty} C_m(\theta; \mle_1) =1$ 
and $\lim_{\theta \rightarrow 0} C_m(\theta; \mle_1) = 0$,
so that the right-side p-value $P_\theta (\mle(Y_1) \geq \mle_1)$ 
leads to a proper distribution function. 
Corresponding confidence density is given by
\begin{align*}
c_{m1}(\theta; \mle_1)
= \frac{1+\sqrt{5}}{2} \frac{\mle_1}{\theta^2} \phi \left( \frac{1+\sqrt{5}}{2} \frac{\mle_1}{\theta} -1 \right)
- \frac{1-\sqrt{5}}{2} \frac{\mle_1}{\theta^2} \phi \left( \frac{1-\sqrt{5}}{2} \frac{\mle_1}{\theta} -1 \right).
\end{align*}
The implied prior based on $\that_1$ is $$
c_{0mi}(\theta; \mle_i) \propto c_{mi}(\theta; \mle_i) / L(\theta; \mle_i) \propto \theta^{-1}
$$
where $L(\theta; \mle_i)$ is the likelihood function from $p_\theta (\mle_i)$. 

\begin{figure}[h!]
\centering
\includegraphics[width=15cm]{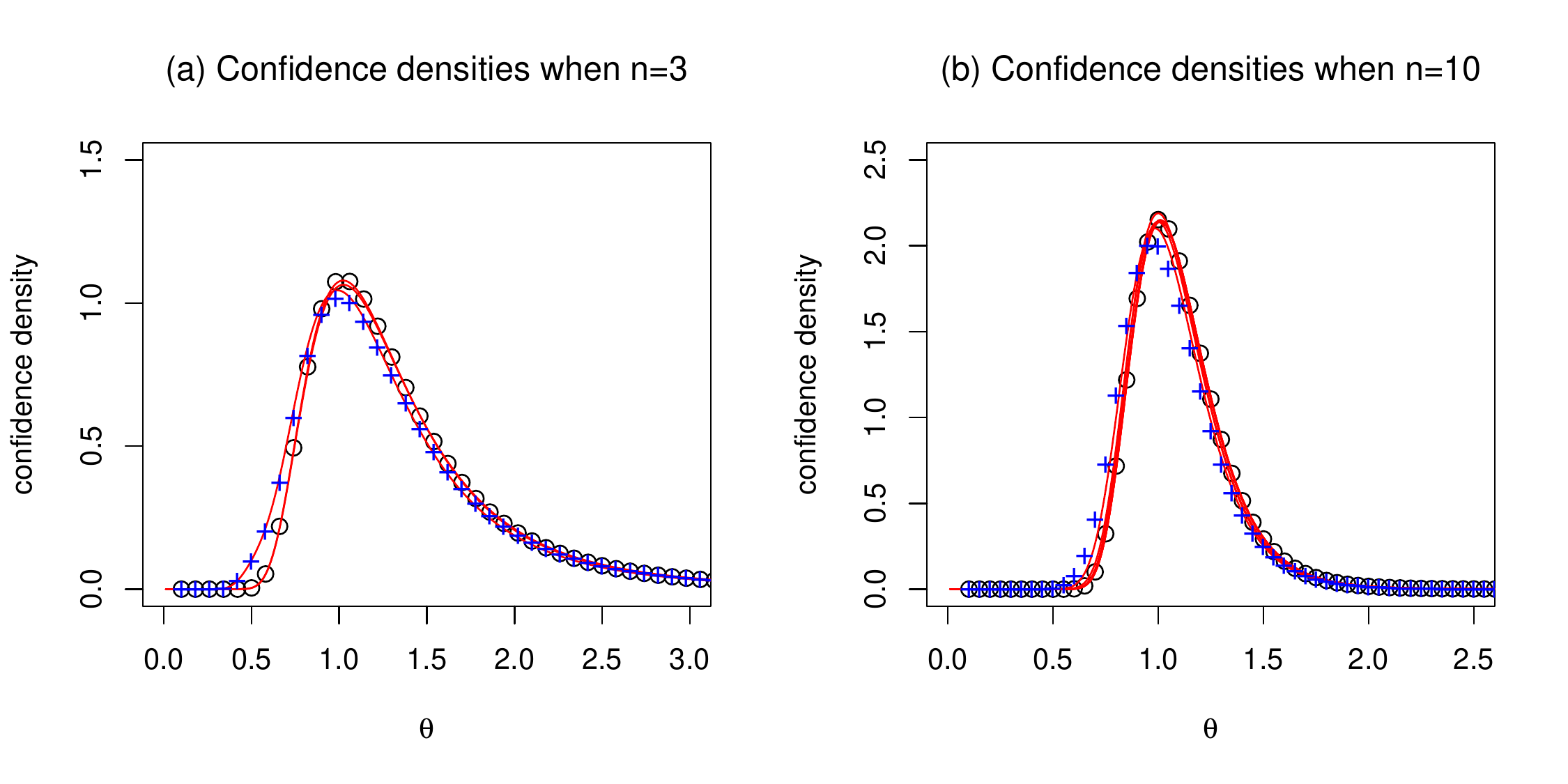}  
\caption{Confidence densities $c_{fi}(\theta; \mle_i)$ for $i=1,2,\cdots, n$ (solid), 
$c_f(\theta; y)$ (circle), $c_m(\theta; \mle_1, \cdots, \mle_n)$ (cross).
(a) $n=3$, (b) $n=10$.}
\label{mu2_conf}
\end{figure}

On the other hand, if we construct the full confidence densities by
$$
c_{fi}(\theta; \mle(y_i), y_{(-i)}) \propto c_{mi}(\theta; \mle(y_i)) L(\theta; y_{(-i)}),
$$
then the resulting confidence density depends on the choice of $y_i$. In this case we should consider $c_{fi}(\theta; \mle_i, y_{(-i)})$ as an approximation to $c_f(\theta;y )$.
Figure~\ref{mu2_conf} plots $n$ confidence densities $c_{fi}(\theta; \mle_i, y_{(-i)})$ (solid) 
and $c_f(\theta; y)$ (circle) with $y_i$ from $N(1,1)$. 
As shown in (b), when $n$ becomes large,
the difference becomes negligible 
and $c_{fi}(\theta; \mle_i, y_{(-i)})$ gets closer to $c_f(\theta; y)$ (circle).

So
$$
c_{fi}(\theta; \mle(y_i), y_{(-i)}) \propto \theta^{-1} L(\theta; \mle_i) L(\theta; y_{(-i)}) 
= \frac{L(\theta; \mle_i)}{L(\theta; y_i)} c_0(\theta)L(\theta;y).
$$
There is a loss of information caused by using $c_{mi}(\theta)$, due to the sign of $y_i$ as captured by $L(\theta; a_i)$. This is negligible even in small samples; see Figure~\ref{mu2_conf}.
However, the marginal confidence
$$c_m(\theta; \mle_1, \cdots, \mle_n) \propto \theta^{-1}L(\theta; \mle_1, \cdots, \mle_n)$$
has a larger loss of information, as shown in both Figure~\ref{mu2_conf}(a) and (b).

\begin{figure}[h!]
\centering
\includegraphics[width=11cm]{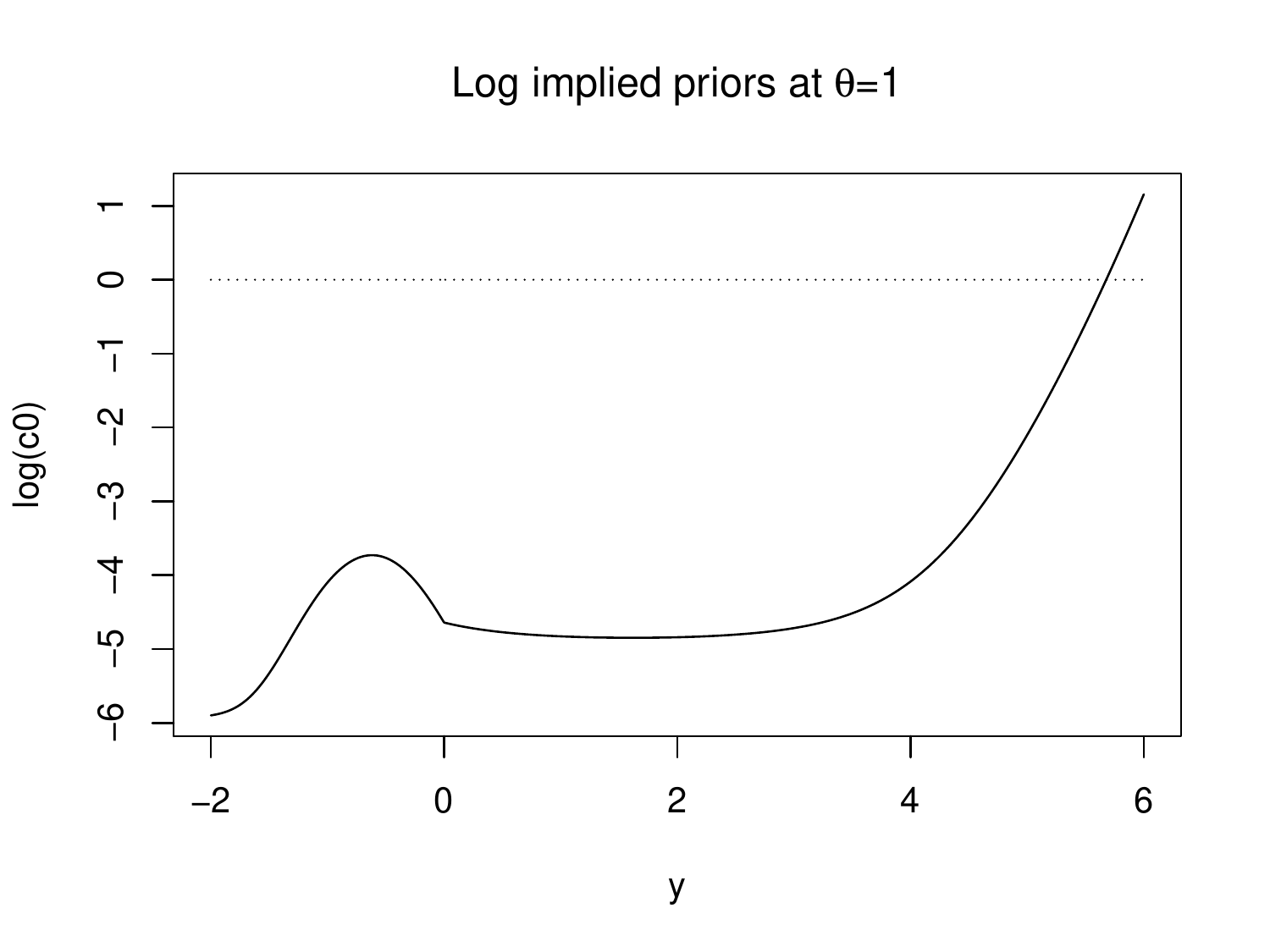}  
\caption{$\log(c_{0m1}(\theta; \mle_1))$ (dotted) and $\log(q_1(\theta; y_1))$ (solid) against $y_1$ varying.}
\label{mu2_logc0}
\end{figure}

Figure~\ref{mu2_logc0} plots the logarithms of
implied prior $c_{0m1}(\theta; \mle_1) \propto 1/\theta$ (dotted)
and $q_1(\theta; y_1) \propto c_{m1}(\theta; \mle_1) / L(\theta; y_1)$ (solid), properly scaled.
Note that $q_1$ is not uniformly continuous on $y_1$,
because the information in $L(\theta; y_1)\propto L(\theta; \mle_1|a_1)$ and $L(\theta; \mle_1)$ differ.

It is also possible to compute the conditional confidence density by using Barndorff-Nielsen's formula (\ref{eq:p-formula}), and to show that we end up with the same implied prior $c_0(\theta)=1/\theta$. Firstly, the likelihood ratio is given by
\begin{align*}
\frac{L(\theta)}{L(\mle)} = \frac{\mle^n}{\theta^n}
\exp \left[ -\frac{b+2n}{4} \left( \frac{\mle^2}{\theta^2} -1 \right)
+\frac{b}{2} \left( \frac{\mle}{\theta} -1 \right)
\right],
\end{align*}
where $b\equiv a^2+a\sqrt{a^2+4n}$, 
and the observed Fisher information 
\begin{align*}
I(\mle) = - \frac{\partial^2 \log L(\theta)}{\partial \theta^2}\Big|_{\theta=\mle}
=\frac{b+4n}{2\mle^2}.
\end{align*}
Then we have
\begin{align*}
p_\theta (\mle | a) \approx & \frac{c}{\sqrt{2}} \frac{\sqrt{b+4n}}{\mle} 
\frac{\mle^n}{\theta^n}
\exp \left[ -\frac{b+2n}{4} \left( \frac{\mle^2}{\theta^2} -1 \right)
+\frac{b}{2} \left( \frac{\mle}{\theta} -1 \right)
\right].
\end{align*}
Let $U \equiv \mle(Y)/\theta$ and let $u=\mle(y)/\theta$, then the conditional density of $u|a$ becomes
\begin{align*}
p_\theta (u | a) \approx & \frac{c\sqrt{b+4n}}{\sqrt{2}} u^{n-1}
\exp \left[ -\frac{b+2n}{4} \left( u^2 -1 \right)
+\frac{b}{2} \left( u-1 \right)
\right] ,
\end{align*}
which does not contain $\theta$. Let $F_a (u) = \int p(u|a)du$, then
\begin{align*}
C_c(\theta; \mle|a) = P_\theta(\mle(Y)\geq \mle) = P_\theta (U \geq \mle/\theta)
=1-F_a(\mle/\theta).
\end{align*}
It gives the conditional confidence density
\begin{align*}
c_c(\theta; \mle|a) = -\frac{\partial F_a(\mle/\theta)}{\partial \theta}
\approx \frac{c\sqrt{b+4n}}{\sqrt{2}} \frac{\mle^n}{\theta^{n+1}}
\exp \left[ -\frac{b+2n}{4} \left( \frac{\mle^2}{\theta^2} -1 \right)F
+\frac{b}{2} \left( \frac{\mle}{\theta} -1 \right)
\right],
\end{align*}
and implied prior $c_0(\theta; \mle|a) \propto c_c(\theta; \mle|a) / L(\theta; y) \propto 1/\theta$.
Thus, the conditional confidence density from Barndorff-Nielsen's formula becomes 
 $$c_c(\theta; \mle|a) = c_f(\theta; y) \propto \theta^{-1}L(\theta;y),$$
which is the same as the full confidence \eqref{eq:mu2_cf}.

\subsection*{A3. Discrete case}

A complication arises in the discrete case since the definition of the P-value is not unique, 
and the coverage probability function is guaranteed not to match any chosen confidence level. 
Given the observed statistic $T=t$, among several candidates, the mid P-value
\begin{equation*}
P_{\theta}(T> t)+\frac{1}{2} P_{\theta}(T=t)
\end{equation*}
is often considered the most appropriate (Lancaster, 1961). 

We shall discuss the specific case of the binomial and negative binomial models: $Y_{1}\sim Bin(n,\theta )$ and $Y_{2}\sim NB(y,\theta )$.
The two models have an identical likelihood, proportional to $L(\theta )\propto
\theta ^{y}(1-\theta )^{n-y},$ but have different probability mass
functions, respectively 
\begin{equation*}
P_{\theta }(Y_{1}=y)={\binom{n}{y}}\theta ^{y}(1-\theta )^{n-y}\text{ and }%
P_{\theta }(Y_{2}=n)={\binom{n-1}{y-1}}\theta ^{y}(1-\theta )^{n-y} 
\end{equation*}
Thus, they have the common MLE $\widehat{\theta }=\widehat{\theta }_{1}=%
\widehat{\theta }_{2}=y/n$. However, the two MLEs have different supports 
\begin{equation*}
\widehat{\theta }_{1}\in \left\{ 0,\frac{1}{n},\cdots ,\frac{n-1}{n}%
,1\right\} \text{ and }\widehat{\theta }_{2}\in \left\{ 1,\frac{y}{y+1},%
\frac{y}{y+2}\cdots \right\} , 
\end{equation*}
and therefore $\widehat{\theta }_{1}$ and $\widehat{\theta }_{2}$ have
different distribution, which lead to different P-values. 
Statistical models such as the binomial and negative
binomial models describe how the unobserved future data will be generated. Thus, all the information about $\theta $ in the data and in the statistical
model is in the extended likelihood. The use of the mid P-values 
\begin{align*}
C(\theta ;y_{1}=y) &=\frac{1}{2}P_{\theta }\left( \widehat{\theta }_{1}=%
\frac{y}{n}\right) +P_{\theta }\left( \widehat{\theta }_{1}>\frac{y}{n}%
\right) =\frac{1}{2} \Big( I_\theta(y, n-y+1) + I_\theta(y+1, n-y) \Big) \\
C(\theta ;y_{2}=n) &=\frac{1}{2}P_{\theta }\left( \widehat{\theta }_{2}=%
\frac{y}{n}\right) +P_{\theta }\left( \widehat{\theta }_{2}>\frac{y}{n}%
\right) =\frac{1}{2} \Big( I_\theta(y, n-y+1) + I_\theta(y, n-y) \Big)
\end{align*}
lead to different confidence densities 
\begin{align*}
c(\theta ;y_{1}) =& \frac{1}{2} \left( \frac{\theta^{y-1}(1-\theta)^{n-y}}{%
B(y, n-y+1)} + \frac{\theta^{y}(1-\theta)^{n-y-1}}{B(y+1, n-y)} \right) \\
c(\theta ;y_{2})=& \frac{1}{2} \left( \frac{\theta^{y-1}(1-\theta)^{n-y}}{%
B(y, n-y+1)} + \frac{\theta^{y-1}(1-\theta)^{n-y-1}}{B(y, n-y)} \right)
\end{align*}
where $I_\theta(\cdot,\cdot)$ is the regularized incomplete beta function and $B(\cdot, \cdot)$ is the beta function.

\begin{figure}[h]
\centering
\includegraphics[width=15cm]{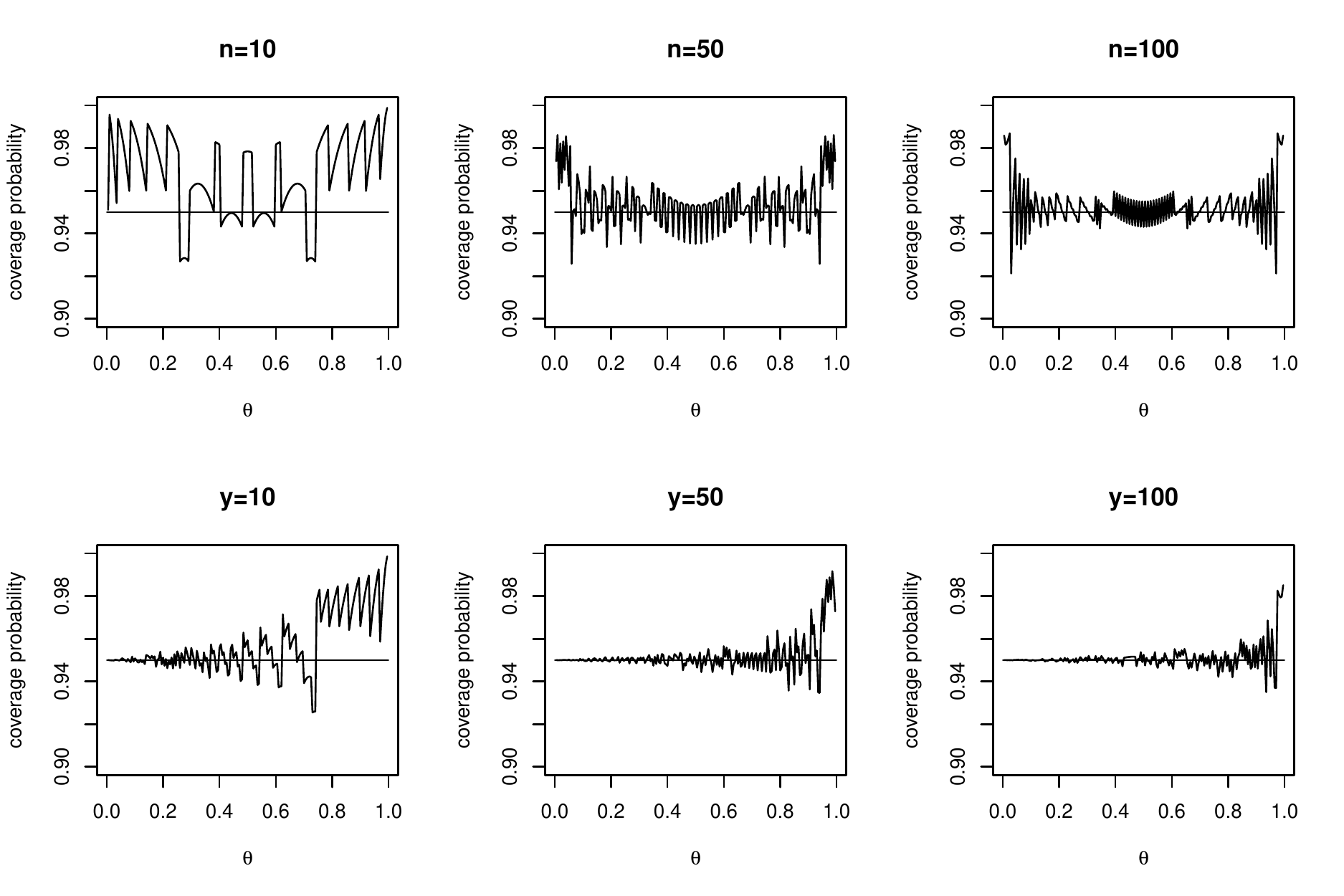}  
\caption{Coverage probabilities of the intervals from the confidence distribution 
based on the mid p-value for binomial models at $n=10, 50, 100$ (top)
and negative binomial models at $y=10, 50, 100$ (bottom).}
\label{fig_add}
\end{figure}

Figure~\ref{fig_add} shows the coverage probabilities of the 95\% two-sided confidence procedure
based on the mid p-value of $\mle$ for binomial models and negative binomial models.
We can see that the coverage probabilities fluctuate around 0.95
but they are not consistently biased in one direction.
Moreover, as $n$ or $y$ becomes larger, 
the difference between the coverage probability and the confidence becomes smaller.
In discrete case, it is not possible to access the exact objective coverage probability of the CI procedure.
Here the confidence is a consistent estimate of the objective coverage probability.
In negative binomial models with $\theta=1$, $y=n$ with probability 1, 
so that it behaves like binomial confidence procedure for $n=y$.

Besides information in the likelihood, the confidence uses
information from the statistical model. 
Consider two different
statistical models, M1: $N=X+1$ where $X \sim$ Poisson($\eta _{1})$ and $%
Y_{1}|N=n\sim Bin(n,\theta )$ and M2: $Y=X+1$ where $X \sim$ Poisson($\eta
_{2})$ and $Y_{2}|Y=y\sim NB(y,\theta ).$ In M1, $Y_{1}|N=n$ and in M2, $%
Y_{2}|Y=y$ have common likelihood, but they are different models, so
that they have no reason to have a common confidence.

\subsection*{A4. When maximal ancillaries are not unique}

When maximal ancillary is not unique, the conditional coverage probability may depend on the choice of the ancillary. However, the lack of unique ancillary does not affect the validity of  Corollary 3 on the absence of relevant subset. We illustrate here with an example from Evans
(2013). The data $y=(y_1,y_2)$ are sampled from a joint distribution with
probabilities under $\theta$ given in the following table:
\begin{center}
\begin{tabular}{ccccc}
$(y_1,y_2)$ & $(1,1)$ & $(1,2)$ & $(2,1)$ & $(2,2)$ \\ \hline
$\theta=1$ & 1/6 & 1/6 & 2/6 & 2/6 \\
$\theta=2$ & 1/12 & 3/12 & 5/12 & 3/12 \\ \hline
\end{tabular}
\end{center}

Here both the data $y$ and parameter $\theta $ are discrete. Strictly, our theory does not cover this case, but we shall use it because it can still illustrate clearly the issues with non-unique maximal ancillaries. The
marginal probabilities are
\begin{eqnarray*}
P_\theta(Y_1=1) = 1/3;\ \ P_\theta(Y_1=2) = 2/3 \\
P_\theta(Y_2=1) = P_\theta(Y_2=2) = 1/2,
\end{eqnarray*}
for $\theta=1,2$. So both $Y_1$ and $Y_2$ are ancillaries, i.e., their
probabilities do not depend on $\theta$. The conditional probabilities of $%
(y_1,y_2)$ given $Y_1=1$ are
\begin{center}
\begin{tabular}{ccc}
$(y_1,y_2)$ & $(1,1)$ & $(1,2)$ \\ \hline
$\theta=1$ & 1/2 & 1/2 \\
$\theta=2$ & 1/4 & 3/4 \\ \hline
\end{tabular}
\end{center}
and, given $Y_2=1$ are
\begin{center}
\begin{tabular}{ccc}
$(y_1,y_2)$ & $(1,1)$ & $(2,1)$ \\ \hline
$\theta=1$ & 1/3 & 2/3 \\
$\theta=2$ & 1/6 & 5/6. \\ \hline
\end{tabular}
\end{center}

Based on the unconditional model, on observing $(y_{1},y_{2})=(1,1)$, we have the likelihood function $L(\theta =1)=1$ and $L(\theta =2)=1/2$, so the MLE $\widehat{\theta}=1$. For $(y_{1},y_{2})=(2,2)$ we have
a different likelihood, but still $\widehat{\theta}=1.$ This means we cannot
reconstruct the likelihood based on the MLE alone, hence MLE is not sufficient. But we can see immediately that we get the same likelihood function under the conditional model given $Y_{1}=1$ or given $Y_{2}=1$, so conditioning on each ancillary recovers the full likelihood and each ancillary is maximal.

Now consider using the MLE itself as a `CI'. Conditional on the ancillaries, the probability that the MLE is correct is
\begin{align*}
& P_{1}(\widehat{\theta }=\theta |Y_{1}=1)=P_{1}(Y_{2}=1|Y_{1}=1)=1/2 \\
& P_{2}(\widehat{\theta }=\theta |Y_{1}=1)=P_{2}(Y_{2}=2|Y_{1}=1)=3/4 \\
& P_{1}(\widehat{\theta }=\theta |Y_{2}=1)=P_{1}(Y_{2}=2|Y_{2}=1)=1/3 \\
& P_{2}(\widehat{\theta }=\theta |Y_{2}=1)=P_{2}(Y_{2}=1|Y_{2}=1)=5/6
\end{align*}
These conditional `coverage probabilities' are indeed distinct from each other.  However, comparing the conditional coverage probabilities given $Y_1$ to that given $Y_2$, \emph{there is no consistent non-trivial bias in one direction} across $\theta$.  So if you use $Y_1$ as the ancillary, you cannot construct further relevant subsets based on $Y_2$. This is the essence of our remark after Corollary 3 that the lack of unique maximal ancillary does not affect the validity of the corollary.

Unfortunately, in this example, the P-value is not defined because the parameter $\theta $ can be an unordered label. So it is not possible to compute any version of confidence function or any implied prior. In the continuous case, we define CI to satisfy $P_{\theta }(\theta \in \mbox{CI})=\gamma $ for all $\theta $. However, in discrete cases, it is often not possible for the coverage probabilities to be same for all $\theta$, which violates the condition of Theorem 1.  Fisher (1973, Chapter III) suggested that for problems such as this, the structure is not sufficient to allow an unambiguous probability-based inference, so only the likelihood is available. 
\end{document}